\documentclass[12pt,reqno]{amsart}
\usepackage{amsmath}
\usepackage{amssymb}
\usepackage{amsthm}
\usepackage{mathrsfs}
\usepackage{latexsym}
\usepackage{xcolor}
\usepackage{comment}

\newtheorem{theorem}{Theorem}[section]
\newtheorem{lemma}{Lemma}[section]
\newtheorem{corollary}{Corollary}[section]
\newtheorem{remark}{Remark}[section]
\newtheorem{definition}{Definition}[section]
\newtheorem{proposition}{Proposition}[section]
\newtheorem{example}{Example}[section]
\newtheorem{assumption}{Assumption}[section]
\numberwithin{equation}{section}
\newcommand{\bth}{\begin{theorem}}
\newcommand{\ethe}{\end{theorem}}

\newcommand{\bre}{\begin{remark}}
\newcommand{\ere}{\end{remark}}

\newcommand{\ble}{\begin{lemma}}
\newcommand{\ele}{\end{lemma}}
\newcommand{\bde}{\begin{definition}}
\newcommand{\ede}{\end{definition}}

\newcommand{\bco}{\begin{corollary}}
\newcommand{\eco}{\end{corollary}}

\newcommand{\bpr}{\begin{proposition}}
\newcommand{\epr}{\end{proposition}}

\newcommand{\bexer}{\begin{exercise}}
\newcommand{\eexer}{\end{exercise}}
\newcommand{\breh}{\begin{hint}}
\newcommand{\ereh}{\end{hint}}

\newcommand{\halmos}{\hfill \qed}

\newcommand{\bexam}{\begin{example}}
\newcommand{\eexam}{\end{example}}

\newcommand{\pr} {{\bf Proof.}}

\newcommand{\bfi}{\begin{fig}}
\newcommand{\efi}{\end{fig}}

\newcommand{\beao}{\begin{eqnarray*}}
\newcommand{\eeao}{\end{eqnarray*}\noindent}

\newcommand{\beam}{\begin{eqnarray}}
\newcommand{\eeam}{\end{eqnarray}\noindent}
\newcommand{\E}{\mathbf{E}}
\newcommand{\PP}{\mathbf{P}}

\newcommand{\xto}{x\to\infty}

\newcommand{\uto}{u\to\infty}

\newcommand{\bG}{\overline{G}}
\newcommand{\bV}{\overline{V}}

\newcommand{\bbr}{{\mathbb R}}

\newcommand{\bbb}{{\mathbb B}}

\newcommand{\bbn}{{\mathbb N}}

\newcommand{\vep}{\varepsilon}

\allowdisplaybreaks[1]
\begin{document}
\title[Insurance and financial risks in a non L\'{e}vy-renewal environment]{Interplay of insurance and financial risks in a non 
L\'{e}vy-Renewal environment}

\author[D.G. Konstantinides, C.D. Passalidis]{Dimitrios G. Konstantinides, Charalampos  D. Passalidis} 

\address{Dept. of Statistics and Actuarial-Financial Mathematics,
University of the Aegean, Karlovassi, GR-83 200 Samos, Greece}
\email{konstant@aegean.gr}
\email{sasd24009@sas.aegean.gr.}

\date{{\small \today}}

\begin{abstract}
In this paper we consider a multivariate risk model, with common counting process and common process of logarithmic returns for the investment portfolio. We assume that the claim-vectors, the counting process and the logarithmic returns of the investment portfolio satisfy a weak dependence structure. Further, we consider that the counting process represents an inhomogeneous renewal process, and the logarithmic returns represent a c\'{a}dl\'{a}g process with independent but not necessarily stationary increments. Under these conditions we provide an asymptotic expression for the infinite-time entrance probability of the discounted aggregate claims into some rare set $x\,A$, where $A$ denotes a set from a general set family, crucial for the actuarial practice, when the common distribution of the claim vectors belong to a multivariate heavy-tailed distribution class, denoted by $(\mathcal{D} \cap \mathcal{A})_A$. This result, is derived under a moment condition for the financial risks, and underlines the multivariate linear single big jump principle. When we restrict the distribution class of the claim-vectors to multivariate regular variation, we find more explicit asymptotic expressions, weakening the moment conditions on the financial risks. The asymptotic formulas, derived through double 'dependence solution', become more direct and practical in applications. With respect to the technical part, due to non L\'{e}vy-Renewal framework, the classical  Kesten-Goldie theorems are not applicable, nor their extensions. The way we make the discretization of the process of the discounted aggregate claims permits to derive uniform asymptotics with respect to the number of summands, that facilitate the approximation of the infinite sums of the main results.
\end{abstract}

\maketitle
\textit{Keywords: Multivariate risk model; Dependent insurance and financial risks; 
Heavy tailed random vectors; Infinite time horizon; Uniformity}
\vspace{3mm}

\textit{Mathematics Subject Classification}: Primary 62P05 ;\quad Secondary 60G70.


\section{Introduction and motivation} \label{sec.KLP.1}

In this paper we consider an insurer, who operates $d$-lines of business, with $d \in \bbn$, that share a common counting process. Namely, the claim vectors of the $d$-lines of business, are described by a sequence of non-negative random vectors $\{{\bf X}^{(i)}\,,\;i \in \bbn\}$, which arrive at the corresponding moments $\{\tau_{i}\,,\;i \in \bbn\}$, with $\tau_0=0$. We notice that to avoid trivial cases, each  vector ${\bf X}^{(i)}=(X_1^{(i)},\,\ldots,\,X_d^{(i)})^{\top}$, can contain zero components, but not all of the components to be zero.

Further, the arrival times  $\{\tau_{i}\,,\;i \in \bbn\}$ represent a counting process  $\{N(t)\,,\;t \geq 0\}$, which is defined as 
\beao
N(t) = \sup \{i \in \bbn\;:\;\tau_i \leq t\}\,,
\eeao 
for any $t \geq 0$, with $\sup \emptyset =0$ conventionally, and has finite mean
\beao
\lambda(t) = \E[N(t)] =\sum_{i=1}^{\infty} \PP(\tau_i \leq t)\,,
\eeao
for any fixed $t \geq 0$. We also define the interval $\Lambda = \{t\;:\;\PP(\theta_1 \leq t) > 0\}$, where the sequence $\{\theta_i := \tau_i - \tau_{i-1}\,,\;i \in \bbn\}$ represents the inter-arrival times between the successive claim vector arrivals. We also assume that the insurer invests his surplus into risk-free or risky investments, or both. The logarithmic returns of the investment portfolio are depicted by a c\'{a}dl\'{a}g stochastic process $\{\xi(t)\,,\;t \geq 0\}$, with initial value $\xi(0)=0$. Then, the discounted aggregate claims up to time $t\geq 0$, are given through the relation
\beao
{\bf D}(t) =\sum_{i=1}^{N(t)} {\bf X}^{(i)}\,e^{-\xi(\tau_i)}=\left( 
\sum_{i=1}^{N(t)} X_{1}^{(i)}\,e^{-\xi(\tau_{i}) } ,\, 
\ldots ,\,
\sum_{i=1}^{N(t)} X_{d}^{(i)}\,e^{-\xi(\tau_{i}) } 
\right)^{\top}\,,
\eeao
for any $t\geq 0$  (where ${\bf y}^T$ denotes the transpose of vector ${\bf y}$), while the corresponding discounted aggregate claims over infinite time horizon, are given by the relation
\beam \label{eq.KLP.1.2}
{\bf D}(\infty) =\sum_{i=1}^{\infty} {\bf X}^{(i)}\,e^{-\xi(\tau_i)}=\left( 
\sum_{i=1}^{\infty} X_{1}^{(i)}\,e^{-\xi(\tau_{i}) } ,\, 
\ldots ,\, 
\sum_{i=1}^{\infty} X_{d}^{(i)}\,e^{-\xi(\tau_{i}) } 
\right)^{\top}\,.
\eeam
We focus our attention on the asymptotic behavior of the probability
\beam \label{eq.KLP.1.3} 
\PP\left({\bf D}(\infty) \in x\,A \right)\,,
\eeam
as $\xto$, where $A$ denotes a set from a general family of sets, and it can take several interesting forms for the actuarial practice (see the set family $\mathscr{R}$, and related discussions in subsection 2.1). This set is immediately connected with the multivariate distribution classes, which is assumed for the common distribution $F$ of the claim vectors.

Before proceeding to more details about the assumptions of our model, it is expedient to overview the recent literature about the estimation of the probability \eqref{eq.KLP.1.3}, and of the corresponding asymptotic relations for ${\bf D}(t)$, with $t< \infty$, or of the asymptotic relations of the ruin probabilities, under different forms of $A$, to show the novelty of this paper. 

The papers \cite{guo:wang:2013}, \cite{cheng:konstantinides:wang:2022}, \cite{chen:li:cheng:2023}, \cite{li:2023}, \cite{yang:chen:yuen:2024}, \cite{chen:konstantinides:passalidis:2025} and \cite{konstantinides:passalidis:2024j}  are only some of the contributions on this topic, in one-dimensional or multidimensional set up, mostly with concrete form of set $A$.

However, in the previous papers consider independent insurance and financial risks, namely the sequences  $\{{\bf X}^{(i)}\,,\;i \in \bbn\}$ and $\{\xi(t)\,,\;t \geq 0\}$ are independent. It is known that these two fundamental risks for the modern insurance industry, usually possess some dependence structure, hence the extension to study risk models with dependent insurance and financial risks is not only of theoretical value (generalization of independent results or counter examples), but also provides important tools for the actuarial practice.

In discrete time risk models, dependence between insurance and financial risks was studied during the last decade, see for example \cite{zhou:wang:wang:2012}, \cite{yang:wang:2013}, \cite{yang:gao:li:2016}, \cite{chen:2017}, \cite{chen:yuan:2017}, \cite{tang:yang:2019} among others.
 
However, in continuous time risk models, we know only five contributions to this topic. For finite time horizon, and one-dimensional set up, \cite{guo:2022} established asymptotic estimation for the ruin probability, in a Poisson risk model, in which the insurance and financial risk processes are two jump diffusion processes and the dependence between the two risks stems from the dependence between the claims and the jumps of the investment portfolio. In \cite{yang:fan:yuen:2023} we find extension of the previous one in renewal risk model, and some other distributional generalizations. In similar line we find  in \cite{cheng:wang:2025}, an application in a reinsurance scheme. We should note that such kind of models have well fit to description of situations like economical crises or pandemics, but do not work in periods of economic stability, in which works our approach here, but only for infinite time horizon.

On multidimensional set up, over infinite time horizon, \cite{cheng:konstantinides:wang:2024} considered a risk model, in which the counting process $\{N(t)\,,\;t \geq 0\}$ represents a renewal process, and the $\{\xi(t)\,,\;t \geq 0\}$ is a multivariate L\'{e}vy process. In that model was assumed that at each renewal 
epoch the distribution of the product of insurance and financial risks belongs to the class of multivariate regular variation, and the components of the product are asymptotically dependent. However, the dependence between the two risks is arbitrary. In \cite{konstantinides:passalidis:2025o}, following the same direction, but considering a common L\'{e}vy process $\{\xi(t)\,,\;t \geq 0\}$ in the $d$-lines of business, the condition of asymptotic dependence was relaxed, and was used distributions from more general multivariate classes. Although there is a variety of assumptions in the previous models, that permit several generalizations, in all these
five paper we find the following two conditions
\begin{enumerate}
\item
We have L\'{e}vy-Renewal risk models, namely the $\{\xi(t)\,,\;t \geq 0\}$ represents a L\'{e}vy process and the  $\{N(t)\,,\;t \geq 0\}$ is a renewal one.
\item
The  $\{N(t)\,,\;t \geq 0\}$ is independent of all the other sources of randomness.
\end{enumerate}   
 
In this paper, we use 
\begin{enumerate}
\item
$\{\xi(t)\,,\;t \geq 0\}$ as a c\'{a}dl\'{a}g process, with independent, but not necessarily stationary increments, while the $\{\theta_{i}\,,\;i \in \bbn\}$ are independent but not necessarily identically distributed, non-negative random variables (hence $\{N(t)\,,\; t\geq 0\}$ is a inhomogeneous renewal process). 
\item
a general enough dependence structure, to describe the interdependence of increments among the sequences $\{{\bf X}^{(i)}\,,\;i \in \bbn\}$, $\{\xi(t)\,,\;t \geq 0\}$ and  $\{N(t)\,,\;t \geq 0\}$,  which contains the independence of any two of the sequences, and of any three of the sequences, as special case.
\end{enumerate}  

The rest of the paper is organized as follows. In Section 2, we provide the preliminaries about the distribution classes and the dependence structures, which used in the paper. In Section 3, we  present the main result together with some corollaries. Further, in Section 4 we give their proofs after some necessary preliminary lemmas. Our techniques are based on the single big jump principle, through discretization of the process of discounted aggregate claims. The way we follow in the proof of our results (see, Lemmas \ref{lem.KLP.3.5}, \ref{lem.KLP.3.6}) permit some kind of uniformity with respect to the number of summands, that plays crucial role in the approximation of asymptotic formulas through Monte Carlo simulations, from the aspect of practitioners.

\section{Distribution classes and dependence structures} \label{sec.KLP.2}

For two real numbers $a,\,b$, we define $a\vee b := \max\{a,\,b\}$, $a\wedge b :=\min\{a,\,b\}$. All the random vectors are of $d$-dimension, denoted through bold script, and their operations are understood component-wisely, namely ${\bf x}+{\bf y}= (x_1 \pm y_1,\,\ldots,\,x_d \pm y_d)^{\top}$, or for a finite quantity $k>0$, the scalar product becomes $k\,{\bf x}=(k\,x_1,\,\ldots,\,k\,x_d)^{\top}$. For any set $\bbb \in \bbr^d:=(-\infty,\,\infty)^d$, we write as $\bbb^c$ for its complement set, as $\overline{\bbb}$ for its closed hull, as $\partial \bbb$ its border, and as ${\bf 1}_{\bbb}$ for its indicator function. 

For any random variable (or, vector) $Y$, we write $Y \stackrel{d}{\sim} V$, when $Y$ follows distribution $V$, and in this case we depict by $s(V)$ or by $s(Y)$ the support of the distribution $V$. If distribution $V$ is one-dimensional, we denote by $\bV(x)=1-V(x)$, for any $x \in \bbr$, its distribution tail. 

Hereafter, all the limit relations hold as $\xto$, except otherwise stated. For two positive, uni-variate functions $f(\cdot)$, $g(\cdot)$, we denote by $f(x) \sim c\,g(x)$, for some $0< c < \infty$, if 
\beao
\lim \dfrac {f(x)}{g(x)} =c\,,
\eeao
by $f(x)=o[g(x)]$, if
\beao
\lim \dfrac {f(x)}{g(x)} =0\,,
\eeao
by $f(x)=O[g(x)]$, if 
\beao
\limsup \dfrac {f(x)}{g(x)} < \infty \,,
\eeao
and $f(x) \asymp g(x)$, if both $f(x)=O[g(x)]$ and $g(x)=O[f(x)]$ are true.
The previous asymptotic notation remains intact for multidimensional set up, for example, for ${\bf f}$, ${\bf g}$ two $d$-variate, positive functions, and $\bbb \in \bbr_+^d$, with ${\bf 0} \notin \overline{\bbb}$. We write $f(x\,\bbb)=O[g(x\,\bbb)]$, if 
\beao
\limsup \dfrac {f(x\,\bbb)}{g(x\,\bbb)} < \infty \,.
\eeao

\subsection{Multivariate heavy tailed distributions}  \label{sec.KLP.2.1}

For sake of compactness of the text, the following definitions are provided for distributions with support on the non-negative quadrant, while all the one-dimensional distributions, have infinite right endpoint (namely it holds $\bV(x) >0$, or $\bG_A(x)>0$, for any $x \in \bbr$).

To define the most of the following multivariate distribution classes, we need the following wide set family
\beam \label{eq.KLP.2.1}
\mathscr{R}=\{ A\subsetneq \bbr^d\;:\;A\;\text{open, increasing}\,,\;A^c\;\text{convex}\,,\; {\bf 0} \notin \overline{A} \}\,,
\eeam 
where ${\bf 0}$ is the origin of the axes, and a set $A$ is called increasing if for any ${\bf x} \in A$, ${\bf y} \in \bbr_+^d$, we have ${\bf x} + {\bf y} \in A$.

\bre \label{rem.KLP.2.1}
We observe that the set $\mathscr{R}$ in \eqref{eq.KLP.2.1} represent a cone with respect to positive scalar multiplication, namely if $A \in \mathscr{R}$ and $k>0$, then we find $k\,A \in \mathscr{R}$. Further, some sets from the family $\mathscr{R}$, that are useful in actuarial practice, are as follows
\beao
A_1=\left\{ {\bf y}\;:\;\sum_{j=1}^d l_j\,y_j > b \right\}\,,
\eeao
with $b>0$ and $l_1,\,\ldots,\,l_d \geq 0$, with $l_1 +\cdots + l_d =1$, and 
\beao
A_2=\left\{ {\bf y}\;:\;y_j > b_j\,,\;\exists\;j=1,\,\ldots,\,d \right\}\,.
\eeao
Hence, the probability in \eqref{eq.KLP.1.3}, when $A=A_i$, with $i=1,\,2$, corresponds to some interesting situations about the solvency of the insurer. More specifically, when $A = A_1$, we find that \eqref{eq.KLP.1.3} represents the probability that the sum of discounted aggregate claims of the $d$-lines of business, exceeds a threshold $b\,x$, while when $A = A_2$, we obtain that \eqref{eq.KLP.1.3} represents the probability that in some of the $d$-lines of business, the discounted aggregate claims exceeds a threshold $b_j\,x$. Even more, in the special sub-case of $A_2$ with $d=1$, we have $A_3 = (b,\,\infty)$, for some $b>0$, with the case $A=(1,\,\infty)$ to show the tail distribution of the discounted aggregate claims in one-dimensional set up. For more details about the family $\mathscr{R}$, see \cite[Sec. 4]{samorodnitsky:sun:2016}.
\ere

Based on the set $\mathscr{R}$, in \cite{samorodnitsky:sun:2016} was introduced the class of multivariate subexponential distributions. For this purpose, was proved that for any random vector ${\bf Z} \stackrel{d}{\sim} G$, with support on the non-negative quadrant, the random variable 
\beao
Z_A=\sup \left\{ u\;:\; {\bf Z} \in u\,A \right\}\,,
\eeao
has proper distribution $G_A$, whose tail is presented by the relation
\beam \label{eq.KLP.2.5}
\bG_A(x)=\PP\left({\bf Z} \in x\,A \right)=\PP\left(\sup_{{\bf p} \in I_A} {\bf p}^{\top}\,{\bf Z} > x \right)\,,
\eeam
for any $x>0$, for some index set $I_A \subsetneq \bbr^d$, see further \cite[Lem. 4.5, Lem. 4.3(c)]{samorodnitsky:sun:2016}, for full argumentation.

Let $A \in \mathscr{R}$ be a fixed set. We say that the distribution $G$ belongs to the class of multivariate subexponential distributions on $A$, symbolically $G \in \mathcal{S}_A$, if $G_A \in \mathcal{S}$, namely for any (or, equivalently, for some) integer $n\geq 2$ it holds
\beao
\lim \dfrac{\overline{G^{n*}_A}(x) }{\bG_A(x)}=n\,,
\eeao
where $G^{n*}_A$ represents the $n$-th fold convolution of $G_A$ with itself. 

By similar way, in \cite{konstantinides:passalidis:2024g} were introduced the classes $\mathcal{L}_A$ and $\mathcal{D}_A$. Concretely, we say that the distribution $G$ belongs to the class of multivariate dominatedly varying distributions on $A$, symbolically $G \in \mathcal{D}_A$, if $G_A \in \mathcal{D}$, namely for any (or, equivalently, for some) $b \in (0,\,1)$ it holds
\beao
\limsup \dfrac{\bG_A(b\,x) }{\bG_A(x)}< \infty\,.
\eeao
And similarly, we say that the distribution $G$ belongs to the class of multivariate long tailed distributions on $A$, symbolically $G \in \mathcal{L}_A$, if $G_A \in \mathcal{L}$, namely for any (or, equivalently, for some) $a >0$ it holds
\beao
\lim \dfrac{\bG_A(x-a) }{\bG_A(x)}=1\,.
\eeao  
In this way, we write $G\in (\mathcal{D} \cap \mathcal{L})_A$ if $G_A \in \mathcal{D} \cap \mathcal{L}$, while from the one-dimensional properties of heavy-tailed distributions, we obtain that $\mathcal{S}_A \subsetneq \mathcal{L}_A$ and $(\mathcal{D} \cap \mathcal{L})_A \equiv (\mathcal{D} \cap \mathcal{S})_A$, see \cite[Ch. 2]{leipus:siaulys:konstantinides:2023} for more information about these inclusions.

In \cite{konstantinides:passalidis:2025h}, was introduced the class of multivariate positively decreasing distributions on $A$, symbolically $(\mathcal{P_D})_A$. We say that $G \in (\mathcal{P_D})_A$, if $G_A \in \mathcal{P_D}$, namely for any (or, equivalently, for some) $v >1$ it holds
\beao
\limsup \dfrac{\bG_A(v\,x) }{\bG_A(x)}<1\,,
\eeao   
as also the class of multivariate subexponential and positively decreasing distributions on $A$, symbolically $\mathcal{A}_A$. We say that $G \in \mathcal{A}_A$, if $G_A \in \mathcal{A}= \mathcal{S} \cap\mathcal{P_D}$. The class of multivariate dominatedly varying, subexponential and positively decreasing distributions on $A$, symbolically $(\mathcal{D} \cap \mathcal{A})_A$, introduced also in the same lines. Namely, we say that $G \in (\mathcal{D} \cap \mathcal{A})_A$, if $G_A \in \mathcal{D} \cap \mathcal{A}= \mathcal{D} \cap \mathcal{P_D}\cap \mathcal{S}$. Let us note that $\mathcal{P_D}$ is a general class, and provides a quite mild restriction on classes $\mathcal{S}$ and $\mathcal{D}  \cap \mathcal{L}$, that means that $\mathcal{A}$ and $\mathcal{D} \cap \mathcal{A}$ are practically negligibly smaller than $\mathcal{S}$ and $\mathcal{D} \cap \mathcal{L}$, respectively, see \cite{tang:2006} and \cite{konstantinides:passalidis:2025b} for more properties of class $\mathcal{P_D}$. For the multidimensional set up, similarly the class $(\mathcal{D} \cap \mathcal{A})_A$ is practically negligibly smaller than $(\mathcal{D} \cap \mathcal{L})_A$. However, the $\mathcal{P_D}$ condition is necessary for the asymptotic estimation of  behavior of probability \eqref{eq.KLP.1.3}. In fact, we need  the existence of finite and positive Matuszewska indexes, to show that the infinite randomly weighted sums are tight, see \cite{yi:chen:su:2011}, \cite{tang:yuan:2016} for more details on this topic.

For all distribution classes above, is defined also the multivariate version over all the $\mathscr{R}$,
\beao
\mathcal{B}_{\mathscr{R}}:= \bigcap_{A \in \mathscr{R}}\,\mathcal{B}_A\,,
\eeao
with $\mathcal{B} \in \{\mathcal{D},\,\mathcal{S},\,\mathcal{L},\,\mathcal{P_D},\,\mathcal{D}\cap \mathcal{L},\,\mathcal{A},\,\mathcal{D} \cap \mathcal{A} \}$. For more discussions about the properties of these classes, can be found in \cite{samorodnitsky:sun:2016}, \cite{konstantinides:passalidis:2024g}, \cite{konstantinides:passalidis:2025h}. For examples of them, see in \cite[Sec. 4]{samorodnitsky:sun:2016}, \cite[Sec. 4]{konstantinides:liu:passalidis:2025}.

Now, we recall some definitions of the regularly varying distributions and the Matuszewska indexes. We say that an one-dimensional distribution $V$ belongs to class of regularly varying distributions, symbolically $V \in \mathcal{R}_{-\alpha}$, with $\alpha \in (0,\,\infty)$, if it holds
\beao
\lim \dfrac {\bV(t\,x)}{\bV(x)} = t^{-\alpha}\,,
\eeao
for any $t>0$.

The lower and upper Matuszewska indexes are defined as
\beam \label{eq.KLP.2.6}
J_V^- = -\lim_{\uto}\dfrac {\log \bV^*(u)}{\log u}\,,\qquad J_V^+ = -\lim_{\uto}\dfrac {\log \bV_*(u)}{\log u}\,,
\eeam
with
\beao
\bV_*(u) = \liminf \dfrac {\bV(u\,x)}{\bV(x)}\,, \qquad \bV^*(u) = \limsup \dfrac {\bV(u\,x)}{\bV(x)}\,.
\eeao
The Matuszewska indexes in \eqref{eq.KLP.2.6} satisfy the relation $0 \leq J_V^- \leq J_V^+ \leq \infty$, while they are related to the characterization of some distribution classes. It is well known that $V \in \mathcal{D}$ if and only if $J_V^+ < \infty$, $V \in \mathcal{P_D}$ if and only if $J_V^- >0$, while in case of $V \in \mathcal{R}_{-\alpha}$ we obtain $J_V^-=J_V^+=\alpha$. For more informations about Matuszewska indexes see \cite[Sec. 2.1.2]{bingham:goldie:teugels:1987}.

Finally, we remind the class of multivariate regular variation, symbolically $MRV$, in its standard from. We say that ${\bf Z} \stackrel{d}{\sim} G$ belongs to the class $MRV$, if there exists an one-dimensional distribution $V \in \mathcal{R}_{-\alpha}$, with $\alpha \in (0,\,\infty)$, and a Radon measure $\mu$, non-degenerate to zero, such that
\beao
\lim \dfrac {1}{\bV(x)}\,\PP\left({\bf Z} \in x\,\bbb \right) = \mu(\bbb)\,,
\eeao
for any Borel set $\bbb \in [0,\,\infty]^d$, with ${\bf 0} \notin \overline{\bbb}$, such that $\mu (\partial \bbb) =0$. In this case we write $G \in MRV(\alpha,\,\mu)$. The previous distribution class was introduced in \cite{haan:resnick:1981}, see also \cite{resnick:2007} for various treatments about $MRV$, and \cite{konstantinides:li:2016}, \cite{li:2016}, \cite{li:2022a}, \cite{chen:liu:2024}, \cite{chen:cheng:zheng:2025}, \cite{yang:xu:fan:2026}, \cite{yuan:lu:fu:2025} for applications on risk theory and risk management.

By \cite[Prop. 3.1]{konstantinides:passalidis:2025h}, we find the inclusions:
\beam \label{eq.KLP.2.8}
MRV \subsetneq (\mathcal{D} \cap \mathcal{A})_{\mathscr{R}} \subsetneq (\mathcal{D} \cap \mathcal{L})_{\mathscr{R}} \subsetneq \mathcal{S}_{\mathscr{R}} \subsetneq  \mathcal{L}_{\mathscr{R}}\,,
\eeam
where $MRV$ denotes all distribution of the form $MRV(\alpha,\,\mu)$. Relation \eqref{eq.KLP.2.8} still holds for distribution classes $\mathcal{B}_A$, for any $A \in \mathscr{R}$, with $\mathcal{B} \in \{\mathcal{D}\cap \mathcal{L},\,\mathcal{D} \cap \mathcal{A},\,\mathcal{S},\,\mathcal{L} \}$. The first four classes in \eqref{eq.KLP.2.8}, and sometimes also the $\mathcal{L}_{\mathscr{R}}$, are related with the multivariate linear single big jump principle, see \cite{konstantinides:passalidis:2024g}. The difference of $MRV$ from the $(\mathcal{D} \cap \mathcal{A})_{\mathscr{R}}$ is not trivial as we also find via the examples in \cite[Sec. 4]{konstantinides:liu:passalidis:2025}. However in spite the wide coverage of class $\mathcal{S}_{\mathscr{R}}$ by the $(\mathcal{D} \cap \mathcal{A})_{\mathscr{R}}$, still there are important multivariate subexponential distributions, for the actuarial practice, which do not belong to the $(\mathcal{D} \cap \mathcal{A})_{\mathscr{R}}$, for example the multivariate Gumbel distributions, see \cite[Exam. 4.5]{konstantinides:liu:passalidis:2025} (let remind that in case of Gumbel it holds $J_G^+ = \infty$).

\subsection{Assumptions of the model}  \label{sec.KLP.2.2}

Now we present the main assumptions with respect to risk model from \eqref{eq.KLP.1.2}. The first assumption is related with the increments of the processes $\{N(t)\,,\;t\geq 0\}$ and $\{\xi(t)\,,\; t\geq 0\}$, that are independent but not necessary identically distributed, and thus the renewal and L\'{e}vy processes are not mandatory, respectively.

\begin{assumption} \label{ass.KLP.2.1}
We suppose that the $\{\theta_i=\tau_i - \tau_{i-1}\,,\;i \in \bbn \}$ are independent, positive, random variables. The $\{\xi(t)\,,\; t\geq 0\}$ is a c\'{a}dl\'{a}g process with independent increments, and $\xi(0)=0$.   
\end{assumption}

The following assumption describes a weak dependence structure among  $\{N(t)\,,\;t\geq 0\}$, $\{{\bf X}^{(i)}\,,\;i \in\bbn \}$ and $\{\xi(t)\,,\; t\geq 0\}$, and represents a non-identical extension of \cite[Ass. 2.2]{konstantinides:passalidis:2025o}. In fact, this assumption in its static form and in one-dimensional set up, is inspired by \cite{asimit:jones:2008}, that was extended later to continuous time models, see \cite{li:tang:wu:2010}, \cite{jiang:wang:chen:xu:2015}, \cite{yuan:lu:2023}, among others.

\begin{assumption} \label{ass.KLP.2.2}
Let $A \in \mathscr{R}$ be a fixed set. We suppose that the $\{({\bf X}^{(i)},\, e^{\xi(\tau_{i-1})-\xi(\tau_i) })\,,\;i \in\bbn \}$ are independent, with respect to $i \in \bbn$, with ${\bf X}^{(i)} \stackrel{d}{\sim} F$ and $e^{\xi(\tau_{i-1})-\xi(\tau_i) }\stackrel{d}{\sim} Q_i$, for any $i \in \bbn$. Further, for any $i\in \bbn$, we assume that there exists some measurable function $h_i\;:\;[0,\,\infty) \to (0,\,\infty)$ such that it holds
\beam \label{eq.KLP.2.9}
0 < \inf_{y \in \Delta_i} h_i(y) \,,
\eeam
where $\Delta_i$ denotes a left area of support $s(Q_i)$, and for which it holds
\beam \label{eq.KLP.2.10}
\PP\left({\bf X}^{(i)} \in x\,A\;\big|\; e^{\xi(\tau_{i-1})-\xi(\tau_i)} =y\right) \sim h_i(y)\,\PP \left({\bf X}^{(i)} \in x\,A\right)\,,
\eeam
uniformly for $y \in s(Q_i)$.   
\end{assumption}

Before the analysis of the assumption above, we should make some technical clarifications.

\bre \label{rem.KLP.2.2}
At first, the uniformity of \eqref{eq.KLP.2.10} can be understood as follows
\beao
\lim \sup_{y \in s(Q_i)} \left|\dfrac{\PP\left({\bf X}^{(i)} \in x\,A\;\big|\; e^{\xi(\tau_{i-1})-\xi(\tau_i)} =y\right)}{h_i(y)\,\PP \left({\bf X}^{(i)} \in x\,A\right)}-1 \right| =0\,.
\eeao
Further, the equality to the conditional terms present a usual technique to simplify the writing, and the left member of \eqref{eq.KLP.2.10} is understood as
\beao
\lim_{\vep \to 0} \PP\left({\bf X}^{(i)} \in x\,A\;\big|\; e^{\xi(\tau_{i-1})-\xi(\tau_i)} \in [y-\vep,\;y+\vep]\right)\,,
\eeao
If $y \in E$, for some $E \neq \emptyset$, such that
\beao
\PP\left(e^{\xi(\tau_{i-1})-\xi(\tau_i)} \in dy \right) = 0\,,
\eeao
then the probability at the left member of \eqref{eq.KLP.2.10}, is understood simply as the unconditional one, and hence $h_i(y) =1$, for any $y \in E$. We also note that the function $h_i(\cdot)$ is bounded from above, namely there exists $C_i \in (0,\,\infty)$ such that it holds
\beam \label{eq.KLP.2.A}
h_{i}(y) \leq C_i\,,
\eeam
for all $y \in s(Q_i)$, see \cite[Prop. 2.4]{cui:wang:2025}. Obviously, if $y \notin s(Q_i)$, then it holds $C_i =1$ according to the above.
\ere

\bre \label{rem.KLP.2.3}
The Assumption \ref{ass.KLP.2.2} includes a dependence structure among the insurance risks, the financial risks and the number of claims. We can see from \eqref{eq.KLP.2.10} that if $\{{\bf X}^{(i)}\,,\;i \in \bbn \}$ is independent of $\{\xi(t)\,,\; t\geq 0\}$ and of $\{N(t)\,,\;t\geq 0\}$ (but the $\{\xi(t)\,,\; t\geq 0\}$ not necessarily independent of $\{N(t)\,,\;t\geq 0\}$), then the \eqref{eq.KLP.2.10} holds with $h(y)=1$, for any $y \in (0,\,\infty)$. However, the opposite is not true, namely if $h(y)=1$, for any $y \in (0,\,\infty)$, this does not implies that the  $\{{\bf X}^{(i)}\,,\;i \in \bbn \}$ is independent of $\{\xi(t)\,,\; t\geq 0\}$ and of $\{N(t)\,,\;t\geq 0\}$ (see \cite[Prop. 2.6]{cui:wang:2025} for a counterexample). Having in mind that $X_A^{(i)} = \sup \{ u \;:\; {\bf X}^{(i)} \in u\,A\}$, $i\in \bbn$, the dependence at \eqref{eq.KLP.2.10} is satisfied in static form by various commonly used copulas. We refer the reader to \cite{li:tang:wu:2010}, \cite{tang:yuan:2016}, \cite{cui:wang:2025}, for examples in this static form of \eqref{eq.KLP.2.10}, where the $h_i$ takes concrete forms.

Hence, Assumption \ref{ass.KLP.2.2} provides a weak dependence structure for the increments of the processes $\{{\bf X}^{(i)}\,,\;i \in \bbn \}$, $\{\xi(t)\,,\; t\geq 0\}$ and $\{N(t)\,,\;t\geq 0\}$, which includes the independence as special case.

The assumption that $\{({\bf X}^{(i)}\,,\;e^{\xi(\tau_{i-1})-\xi(\tau_i) })\,,\;i \in \bbn \}$ are independent, implies that the dependence of these three sources of randomness happens only in the time interval $(\tau_{i-1},\,\tau_i]$, which is plausible enough from practical point of view. Finally, if $\{\xi(t)\,,\;t\geq 0\}$, be a L\'{e}vy process, then the  \eqref{eq.KLP.2.10} is reduced to
\beam \label{eq.KLP.2.12}
\PP\left({\bf X}^{(i)} \in x\,A\;\big|\; e^{-\xi(\theta_i)} =y\right) \sim h_i(y)\,\PP \left({\bf X}^{(i)} \in x\,A\right)\,,
\eeam
uniformly for $y \in s(Q_i)$. If additionally we have a renewal $\{N(t)\,,\;t\geq 0\}$, then relation \eqref{eq.KLP.2.12} holds with the same $h_i$, for all $i \in \bbn$ and the $\theta_i$ can be replaced by the $\theta_1$. 
\ere

\bre \label{rem.KLP.2.4}
The dependence structure in relation \eqref{eq.KLP.2.10}, can be 'solved'. This follows from the fact that for any $i \in \bbn$, it holds
\beam \label{eq.KLP.2.13}
\E\left[ h_i\left( e^{\xi(\tau_{i-1})-\xi(\tau_i) }\right) \right] =1\,,
\eeam
Indeed, to see this it is enough to integrate with respect to $\PP\left(e^{\xi(\tau_{i-1})-\xi(\tau_i) } \in dy \right)$ on $s(Q_i)$ the two members of \eqref{eq.KLP.2.10}.
\ere

Let now consider a new stochastic process $\{ \xi_{h_i}(t)\,,\; t\geq 0 \}$, that is independent of all the other sources of randomness, and such that it holds
\beam \label{eq.KLP.2.14}
\PP\left(e^{\xi_{h_i}(\tau_{i-1}^*)-\xi_{h_i}(\tau_i^*)} \in dy \right) =h_i(y)\,\PP\left(e^{\xi(\tau_{i-1})-\xi(\tau_i) } \in dy \right) \,,
\eeam
where $\{\tau_i^*\,,\; i \in \bbn\}$ is a sequence of independent from all the other sources of randomness with $\tau_i^* \stackrel{d}{=} \tau_i$, for any $i \in \bbn$.

Then, by \eqref{eq.KLP.2.13}, we obtain that $\{\xi_{h_i}(t)\,,\; t\geq 0\}$ represents a c\'{a}dl\'{a}g process, (recall also \eqref{eq.KLP.2.9} and \eqref{eq.KLP.2.A}), and furthermore its increments are also independent.

The dependence of Assumption \ref{ass.KLP.2.2}, eventually seems somehow cumbersome (non-explicit) due to the fact that \eqref{eq.KLP.2.10} contains the terms $e^{\xi(\tau_{i-1})-\xi(\tau_i) }$, that include both processes $\{\xi(t)\,,\; t\geq 0 \}$ and $\{N(t) \,,\; t\geq 0\}$ as well. For this reason, we present the following Example, that gives a sufficient condition for validity of \eqref{eq.KLP.2.10}, by a more natural-dynamic way.

\bexam \label{exam.KLP.2.1}
 Let for any $i \in \bbn$ the $h_i \;:\;[0,\,\infty) \to (0,\,\infty)$
be such that \eqref{eq.KLP.2.9} is valid. If it holds
\beam \label{eq.KLP.2.e.1}
\PP\left({\bf X}^{(i)} \in x\,A\;|\;e^{\xi(s_{i-1})-\xi(s_i)} =y\,,\;\tau_{i=1}=s_{i-1}\,,\;\tau_i = s_i \right) \sim h_i(y)\,\PP\left({\bf X}^{(i)} \in x\,A \right)\,,
\eeam
uniformly for $y \in s(Q_i)$ and $s_{i-1} \in \Lambda$, under the convention $s_{i-1} < s_i < \infty$, then \eqref{eq.KLP.2.10} remains true.

Indeed, if $K:=\{s_{i-1} \in \Lambda\,,\;s_{i-1} < s_i < \infty\}$, then for any $i \in \bbn$ we obtain
\beam \label{eq.KLP.3.34} \notag
&&\PP\left({\bf X}^{(i)} \in x\,A\;\big|\; e^{\xi(\tau_{i-1}) - \xi(\tau_i)}=y\right)= \dfrac{\PP\left({\bf X}^{(i)} \in x\,A\,,\; e^{\xi(\tau_{i-1}) - \xi(\tau_i)}=y\right)}{\PP\left( e^{\xi(\tau_{i-1}) - \xi(\tau_i)}=y\right)}\\[2mm] \notag
&&= \dfrac{\int_{K} \PP\left({\bf X}^{(i)} \in x\,A\,,\; e^{\xi(s_{i-1}) - \xi(s_i)}=y\,,\;\tau_{i-1} \in ds_{i-1}\,,\;\tau_{i} \in ds_{i}\right)}{\int_{K} \PP\left( e^{\xi(s_{i-1}) - \xi(s_i)}=y\,,\;\tau_{i-1} \in ds_{i-1}\,,\;\tau_{i} \in ds_{i}\right)}\\[2mm] \notag
&&=  \dfrac{1}{\int_{K} \PP\left( e^{\xi(s_{i-1}) - \xi(s_i)}=y\,,\;\tau_{i-1} \in ds_{i-1}\,,\;\tau_{i} \in ds_{i}\right)}\\[2mm]
&&\qquad\times \int_{K}\PP\left({\bf X}^{(i)} \in x\,A\;\big|\;e^{\xi(s_{i-1}) - \xi(s_i)}=y\,,\;\tau_{i-1} = s_{i-1}\,,\;\tau_{i} =s_{i}\right)\\[2mm] \notag
&& \qquad \times \PP\left( e^{\xi(s_{i-1}) - \xi(s_i)}=y\,,\;\tau_{i-1} \in ds_{i-1}\,,\;\tau_{i} \in ds_{i}\right)\sim h_i(y)\,\PP\left({\bf X}^{(i)} \in x\,A\right)\,,
\eeam 
where at the last step we used Assumption \ref{ass.KLP.2.2}, via the dominated convergence theorem, through the uniformity of \eqref{eq.KLP.2.10}. 
\eexam

\section{Asymptotics for discounted aggregate claims} \label{sec.KLP.3}

Here we present the main results, that are oriented to the estimation of probability in \eqref{eq.KLP.1.3}, as $\xto$, as also their proofs.

\subsection{Main results} \label{sec.KLP.3.1}

The following statement is the main result of the paper. We recall that if $F_A \in \mathcal{D}\cap \mathcal{A} \subsetneq \mathcal{D}\cap \mathcal{P_D}$, hence it holds $0 < J_{F_A}^- \leq  J_{F_A}^+ < \infty$.

\bth \label{th.KLP.3.1}
Let $A \in \mathscr{R}$ be a fixed set, and the discounted aggregate claims be of the \eqref{eq.KLP.1.2}. We suppose that the Assumptions \ref{ass.KLP.2.1} and \ref{ass.KLP.2.2} hold. If $F \in (\mathcal{D}\cap \mathcal{A})_A$ and there exist $0< p_1< J_{F_A}^- \leq  J_{F_A}^+ < p_2 < \infty$, such that 
\beam \label{eq.KLP.3.1}
\E\left[e^{p_1(\xi(\tau_{i-1})-\xi(\tau_i))} \bigvee e^{p_2(\xi(\tau_{i-1})-\xi(\tau_i))}\right] <1\,,
\eeam
for any $i \in \bbn$, then we obtain
\beam \label{eq.KLP.3.2} \notag
&&\PP\left({\bf D}(\infty) \in x\,A \right) \sim \sum_{i=1}^{\infty} \PP\left( {\bf X}^{(i)}\,e^{-\xi(\tau_i)} \in x\,A \right) \\[2mm]
&& \sim \sum_{i=1}^{\infty} \PP\left( {\bf X}^{(i)}\,e^{\xi_{h_i}(\tau_{i-1}^*) - \xi_{h_i}(\tau_{i}^*)}\,e^{-\xi(\tau_{i-1})} \in x\,A \right) \,.
\eeam
\ethe

\bre \label{rem.KLP.3.1}
We notice that if the assumption of Theorem \ref{th.KLP.3.1}, are valid for any  $A \in \mathscr{R}$, that means $F \in (\mathcal{D}\cap \mathcal{A})_\mathscr{R}$, and further Assumption  \ref{ass.KLP.2.2} holds for any  $A \in \mathscr{R}$, and the moment condition \eqref{eq.KLP.3.1}, holds for any  $A \in \mathscr{R}$, and some real numbers $p_1,\,p_2$ for which $0< p_1< J_{F_A}^- \leq  J_{F_A}^+ < p_2 < \infty$, then relation \eqref{eq.KLP.3.2} is true for any  $A \in \mathscr{R}$. The conditions of Theorem \ref{th.KLP.3.1} are general enough, to permit simultaneous dependence of the claim vectors, of the number of claims, and of the logarithmic returns of the investment portfolio. According to the best of our knowledge, such models have not been studied before, even in one-dimensional set up.
\ere

\bre \label{rem.KLP.3.1,5}
Theorem \ref{th.KLP.3.1} has two 'good' properties with respect to its applicability. The first one is that in the last sum of \eqref{eq.KLP.3.2}, 
each probability  contains products of independent factors, which makes its application much easier.

The second one is that equation \eqref{eq.KLP.3.2} holds in fact uniformly for all $n \in \bbn$, see Lemma \ref{lem.KLP.3.5}, and therefore the infinite sum at the last term of \eqref{eq.KLP.3.2} can be approximated by some finite sum with 'sufficiently large $n \in \bbn$'. This is important for the use of monte-carlo simulations for the approximations of probability \eqref{eq.KLP.1.3}. In literature, such kind of problems (of infinite time) were faced via some change probability measure techniques, see for example \cite{asmussen:glynn:2007}, however such approximations have been used only in more conservative risk models without interest force. So, in more complex risk models, the uniformity with respect to all time horizons (or, here, the uniformity with respect to $n \in \bbn$) functions as antidote for the practitioners (see, \cite{tang:2004} for more discussions on uniformity). 
\ere

\bre \label{rem.KLP.3.2}
The products of insurance and financial risks, at each arrival epoch, 
\beao 
{\bf X}^{(i)}\,e^{\xi(\tau_{i-1})-\xi(\tau_i)} \stackrel{d}{\sim} H_i\,,
\eeao 
with $i \in \bbn$, are such that $H_i \in (\mathcal{D}\cap \mathcal{A})_A$. To show this, we just check that $ X_A^{(i)}\,e^{\xi(\tau_{i-1})-\xi(\tau_i)} \stackrel{d}{\sim} H_A^{(i)}$, with 
\beao
X_A^{(i)} = \sup\{u\;:\; X_A^{(i)} \in u\,A\}\,,
\eeao
for any $i \in \bbn$.

Since $X_A^{(i)} \stackrel{d}{\sim} F_A^{(i)} \in \mathcal{D}\cap \mathcal{A}$, by relation \eqref{eq.KLP.3.11} from the proof of Lemma \ref{lem.KLP.3.1}, we obtain that for the random variables
\beao
W_i = e^{\xi(\tau_{i-1})-\xi(\tau_i)}\,,\qquad W_{h,i}=e^{\xi_{h_i}(\tau_{i-1}^*)-\xi_{h_i}(\tau_i^*)}\,,
\eeao
(recall relation \eqref{eq.KLP.2.14}), it holds 
\beam \label{eq.KLP.3.3}
\PP\left(X_A^{(i)}\,e^{\xi(\tau_{i-1})-\xi(\tau_i)} > x \right) \sim \PP\left(X_A^{(i)}\,e^{\xi_{h_i}(\tau_{i-1}^*)-\xi_{h_i}(\tau_i^*)} > x \right)\,.
\eeam
The financial factor $e^{\xi_{h_i}(\tau_{i-1}^*)-\xi_{h_i}(\tau_i^*)}$ is independent of $X_A^{(i)} \stackrel{d}{\sim} F_A^{(i)} \in \mathcal{D}\cap \mathcal{A}$ and 
\beam \label{eq.KLP.3.4} 
&&\E\left[e^{p_2[\xi_{h_i}(\tau_{i-1})-\xi_{h_i}(\tau_i)]} \right]= \int_0^{\infty} y^{p_2} \PP\left(e^{\xi_{h_i}(\tau_{i-1}^*)-\xi_{h_i}(\tau_i^*)}\in dy \right) \\[2mm] \notag 
&& = \int_0^{\infty} y^{p_2}\,{h_i}(y)\,\PP\left(e^{\xi(\tau_{i-1})-\xi(\tau_i)}\in dy \right) \leq C_i\,\E\left[e^{p_2[\xi(\tau_{i-1})-\xi(\tau_i)]} \right] < C_i\,, 
\eeam
where at the next to last step we used that the function ${h_i}$ is bounded from above (see \eqref{eq.KLP.2.A}), while at the last step we took into 
account the moment condition \eqref{eq.KLP.3.1} for the process $\{ \xi(t)\,,\; t\geq 0 \}$. Hence, via relation \eqref{eq.KLP.3.4}, and applying \cite[Th. 2.2(iii), Th. 3.3(ii)]{cline:samorodnitsky:1994} (see also \cite[Cor. 5.2(c)]{leipus:siaulys:konstantinides:2023}), we find that the distribution of the product $X_A^{(i)}\,e^{\xi_{h_i}(\tau_{i-1})-\xi_{h_i}(\tau_i)}$ belongs to the class $\mathcal{D}\cap \mathcal{L}$, while by \cite[Th. 5.1(i)]{konstantinides:passalidis:2025b}, it also belongs to class $\mathcal{P_D}$. Therefore the distribution of this product belongs to $\mathcal{D}\cap \mathcal{A}$ and consequently by \eqref{eq.KLP.3.3} follows that $H_A^{(i)} \in  \mathcal{D}\cap \mathcal{A}$, or equivalently $H_i \in  (\mathcal{D}\cap \mathcal{A})_A$.

This observation, brings the idea that relation \eqref{eq.KLP.3.1} implies a type of multivariate linear single big jump principle, for the asymptotic behavior of the entrance probability of the discounted aggregate claims into set $x\,A$. 
\ere

\bre \label{rem.KLP.3.3}
From Remark \ref{rem.KLP.3.2}, we can make a comparison with the \cite[Th. 4.1]{konstantinides:passalidis:2025o}, where was established relation \eqref{eq.KLP.3.5}, in a L\'{e}vy-Renewal model, in which the insurance and financial risks, at each arrival epoch, are arbitrarily dependent, with ${\bf X}^{(i)}\,e^{\xi(\tau_{i-1})-\xi(\tau_i)} \stackrel{d}{\sim} H \in  (\mathcal{D}\cap \mathcal{A})_A$, and under some assumption on the Laplace exponent of the  L\'{e}vy process (which implies \eqref{eq.KLP.3.1}). Although in that paper there exists arbitrary dependence between the insurance and financial risks, and in some cases the process $\{ \xi(t)\,,\; t\geq 0 \}$ dominates on the $\{{\bf X}^{(i)}\,,\;i \in\bbn \}$, while now they possess a weak dependence structure of Assumption \ref{ass.KLP.2.2}, here we are not restricted to L\'{e}vy-Renewal models, and further in Theorem \ref{th.KLP.3.1} we permit also time-dependent framework, which is not contained in  \cite[Th. 4.1]{konstantinides:passalidis:2025o}. So Theorem \ref{th.KLP.3.1} and Corollary  \ref{cor.KLP.3.1} does not cover neither is covered by \cite[Th. 4.1]{konstantinides:passalidis:2025o}.
We also note that if $\{ \xi(t)\,,\; t\geq 0 \}$ is a  L\'{e}vy process (and the $\{N(t)\,,\; t\geq 0\}$ inhomogeneous renewal), with Laplace exponent $\phi(\cdot)$, such that $\phi(p_2) < 0$, for some $p_2 > J_{F_A}^+$, then \eqref{eq.KLP.3.1} still is valid.
\ere

In the following corollary we obtain a more explicit expression in comparison with \eqref{eq.KLP.3.2}, when the $F$ is restricted in class $MRV$. Further, this restriction permit us to slightly relax the moment condition \eqref{eq.KLP.3.1}.

\bco \label{cor.KLP.3.a}
Let $A \in \mathscr{R}$ some fixed set, and let the discounted aggregate claims of \eqref{eq.KLP.1.2}. Under the conditions of Theorem \ref{th.KLP.3.1}, with the restriction $F \in MRV(\alpha,\,\mu)$, and $\alpha \in (0,\,\infty)$, and instead of \eqref{eq.KLP.3.1}, we assume that the conditions
\beam \label{eq.KLP.3.a}
\E\left[e^{\alpha\,(\xi(\tau_{i-1})-\xi(\tau_i))}\right] <1\,, \qquad \E\left[e^{p_2(\xi(\tau_{i-1})-\xi(\tau_i))}\right] < \infty\,,
\eeam
hold for some $p_2 > \alpha$, for any $i \in \bbn$. Then we obtain
\beam \label{eq.KLP.3.b} 
\PP\left({\bf D}(\infty) \in x\,A \right) \sim \mu(A) \,\bV(x) \,\sum_{i=1}^{\infty} \E\left[ e^{-\alpha\,(\xi_{h_i}(\tau_{i-1}^*) - \xi_{h_i}(\tau_{i}^*)} \right]\,\E\left[ e^{-\alpha\,\xi(\tau_{i-1})} \right]\,.
\eeam
\eco

\bre
Relation \eqref{eq.KLP.3.b} is explicit enough, since the decay-rate of \eqref{eq.KLP.1.3} is determined only through the one-dimensional distribution $V \in \mathcal{R}_{-\alpha}$, with $\alpha \in (0,\,\infty)$, while the rest quantities in the right member of \eqref{eq.KLP.3.b} are positive constants (remind that for any $A \in \mathscr{R}$, we have $\mu(A) \in (0,\,\infty)$, see the proof of \cite[Prop. 4.14]{samorodnitsky:sun:2016}). Next, as also in \eqref{eq.KLP.3.2} (see, Remark \ref{rem.KLP.3.1,5}), relation \eqref{eq.KLP.3.b} in fact holds uniformly with respect to number of summands $n \in \bbn$, see Lemma \ref{lem.KLP.3.6}.
\ere 

In the following corollary we find an integral form of relation \eqref{eq.KLP.3.2}, for the case when the arrival process $\{N(t)\,,\;t\geq 0\}$ is independent of the rest sources of randomness.

\bco \label{cor.KLP.3.1}
\begin{enumerate}
\item[(i)]
Under assumptions of Theorem \ref{th.KLP.3.1}, with the restriction that $\{\theta_i\,,\;i \in\bbn \}$ are independent of all the other sources of randomness, it holds
\beam \label{eq.KLP.3.5} 
&&\PP\left({\bf D}(\infty) \in x\,A \right) \sim \int_0^{\infty} \PP\left({\bf X}\,e^{-\xi(s)} \in x\,A \right) \lambda(ds) \\[2mm] \notag
&& \sim \sum_{i=1}^{\infty} \int_{K} \PP\left({\bf X}^{(i)}\,e^{\xi_{h_i}(s_{i-1})-\xi_{h_i}(s_{i})}\,e^{-\xi(t)} \in x\,A \right) \,\PP\left(\tau_{i-1}^* \in ds_{i-1}\,,\;\tau_{i}^* \in ds_{i} \right)\, \lambda(dt) \,.
\eeam
\item[(ii)]
Under the conditions of part (i), with the conditions of Corollary \ref{cor.KLP.3.a} instead of Theorem \ref{th.KLP.3.1} it holds
\beam \label{eq.KLP.aa} \notag
&&\PP\left({\bf D}(\infty) \in x\,A \right) \sim \mu(A) \,\bV(x) \,\sum_{i=1}^{\infty} \int_{K} \E\left[ e^{-\alpha\,(\xi(s_{i-1}) - \xi(s_{i}))} \right]\,\E\left[ e^{-\alpha\,\xi(t)} \right]\\[2mm]
&&\times \PP(\tau_{i-1}^* \in ds_{i-1}\,,\;\tau_{i}^* \in ds_{i})\,\lambda(dt)\,.
\eeam
\end{enumerate}
\eco

\section{Proofs of main results}  \label{sec.KLP.3.2}

Our technique of proof is based on discretization of process ${\bf D}(\infty)$, and in the presence of multivariate linear single big jump principle. Hence, we need some preliminary lemmas, mostly on the randomly weighted sums on discrete time, and for this we employ the Assumption \ref{ass.KLP.3.1}. In some sense this approach reminds the usual Kesten-Goldie theorems (see \cite{kesten:1973}, \cite{goldie:1991}) on stochastic recurrence equations, but firstly is not restricted to $MRV$, secondly has not independent and identical innovations (that used on L\'{e}vy-Renewal models), while simultaneously contains exact asymptotic expressions. See, \cite{kabanov:legenkiy:promyslov:2026}, \cite{promyslov:2026} for some recent approaches of infinite-time ruin probability through stochastic recurrence equations, in one-dimensional L\'{e}vy-Renewal risk models with independent insurance and financial risks.
Let us note that for $\{{\bf X}^{(i)}\,,\;i \in \bbn \}$, we denote $X_A^{(i)}=\sup\{u\;:\;{\bf X}^{(i)} \in u\,A\} \stackrel{d}{\sim} F_A^{(i)}$, for any $A \in \mathscr{R}$, and $i \in \bbn$.

\begin{assumption} \label{ass.KLP.3.1}
Let $A \in \mathscr{R}$ be a fixed set. We suppose that the $\{({\bf X}^{(i)},\, W_{i}),\;i \in\bbn \}$ are independent random vectors (with respect to $i \in \bbn$), with ${\bf X}^{(i)} \stackrel{d}{\sim} F$ and $W_i \stackrel{d}{\sim} G_i$, for any $i \in \bbn$. Further, for any $i\in \bbn$, we assume that there exists some measurable function $h_i\;:\;[0,\,\infty) \to (0,\,\infty)$ such that it holds
\beam \label{eq.KLP.3.6}
0 < \inf_{y \in E_i} h_i(y) \,,
\eeam
where $E_i$ denotes a left area of the right endpoint for the support $s(G_i)$, and for which it holds
\beam \label{eq.KLP.3.7}
\PP\left({\bf X}^{(i)} \in x\,A\;\big|\; W_i =y\right) \sim h_i(y)\,\PP \left({\bf X}^{(i)} \in x\,A\right)\,,
\eeam
uniformly for $y \in s(G_i)$.   
\end{assumption}

Recall that the sequence $\{{\bf X}^{(i)},\;i \in\bbn \}$ is called tail asymptotic independent on A, symbolically $TAI_A$, if the sequence $\{X_A^{(i)},\;i \in \bbn \}$ is $TAI$, namely for any $i,\,j\,\in \bbn$, with $i\neq j$ it holds
\beao
\lim_{x_i \wedge x_j \to \infty} \PP(X_A^{(i)} >x_i \;\big|\;X_A^{(j)} >x_j)= 0\,,
\eeao
where $TAI$, was introduced by \cite{geluk:tang:2009}. Hereafter, for any $i \in \bbn$ we denote
\beao
\Pi_i = \prod_{j=1}^i W_j\,, 
\eeao
where conventionally $\prod_{j=1}^0= 1$.

We also introduce $W_{h,i}$, with $i \in \bbn$, that represent random variables independent from all the other sources of randomness (and independent each other), with $W_{h,i} \stackrel{d}{\sim} G_{h,i}$, and $G_{h,i}$ given by the relation
\beam \label{eq.KLP.3.12}
G_{h,i}(dy) = h_i(y)\,G_i(dy)\,,
\eeam
which presents a proper distribution, since $\E[h_i(W_i)] = 1$ (as follows by integration of both 
members of \eqref{eq.KLP.3.7}, with respect to $W_i$, over the entire $s(W_i)$). 

In the first lemma we show that under Assumption \ref{ass.KLP.3.1}, the sequence $\{{\bf X}^{(i)}\,\Pi_i,\;i \in\bbn \}$, with $F \in \mathcal{D}_A$, is $TAI_A$, with distributions in class $\mathcal{D}_A$, while with $F \in (\mathcal{D}\cap \mathcal{L})_A$, then the random variables $\{{\bf X}^{(i)}\,\Pi_i,\;i \in\bbn \}$ have distributions again in $(\mathcal{D}\cap \mathcal{L})_A$. Additionally, we prove a type of weak tail equivalence on terms of the sequence.

\ble \label{lem.KLP.3.1}
Let $A \in \mathscr{R}$ be a fixed set. We suppose that the $\{({\bf X}^{(i)},\, W_{i}),\;i \in\bbn \}$ satisfy Assumption \ref{ass.KLP.3.1}, with $F \in \mathcal{D}_A$ and it holds 
\beao
\E\left[W_i^p \right]< \infty\,,
\eeao 
for any $i \in \bbn$, and for some $p > J_{F_A}^+$. Then:
\begin{enumerate}
\item[(i)]
The terms of $\{{\bf X}^{(i)}\,\Pi_i,\;i \in \bbn \}$ have distributions in class $\mathcal{D}_A$, and it holds
\beam \label{eq.KLP.3.8}
\PP\left({\bf X}^{(i)}\,\Pi_i \in x\,A\right) \asymp \PP \left({\bf X}^{(i)} \in x\,A\right)\,,
\eeam
for any $i \in \bbn$.
\item[(ii)]
It holds
\beam \label{eq.KLP.3.B}
\PP\left({\bf X}^{(i)}\,\Pi_i \in x\,A\right) \sim \PP \left({\bf X}^{(i)}\,W_{h,i}\,\Pi_{i-1}  \in x\,A\right)\,,
\eeam
\item[(iii)]
The sequence $\{{\bf X}^{(i)}\,\Pi_i,\;i \in\bbn \}$ is  $TAI_A$. 
\item[(iv)]
If additionally $F \in(\mathcal{D}\cap \mathcal{L})_A $,  the random variables $\{{\bf X}^{(i)}\,\Pi_i,\;i \in\bbn \}$ have distributions in class $(\mathcal{D}\cap \mathcal{L})_A$. 
\end{enumerate}
\ele

\pr~
First we observe that 
\beam \label{eq.KLP.3.9}
\PP\left({\bf X}^{(i)}\,\Pi_i \in x\,A\right) = \PP \left(\sup_{{\bf p} \in I_A}{\bf p}^T\,({\bf X}^{(i)}\,\Pi_i) > x\right) =\PP\left(X_A^{(i)}\,\Pi_i >x \right)\,.
\eeam
\begin{enumerate}
\item[(i)]
At first we shall show relation \eqref{eq.KLP.3.8}. Indeed, from \eqref{eq.KLP.3.9} it is enough 
to show  that 
\beam \label{eq.KLP.3.10}
\PP\left(X_A^{(i)}\,\Pi_i >x \right) \asymp \PP \left(X_A > x\right)\,,
\eeam
for any $i \in \bbn$. From the condition that there exists some $p> J_{F_A}^+$, such that $\E\left[W_i^p \right]< \infty$, 
we find that 
\beam \label{eq.KLP.3.C}
\PP\left(W_i >x \right)=o\left[ \PP \left(X_A > x\right)\right]\,,
\eeam  
and hence, from \cite[Lem. 4.5(i)]{yang:gao:li:2016} (for only one factor in the product), 
we obtain that
\beam \label{eq.KLP.3.11}
\PP\left(X_A^{(i)}\,W_i >x \right) \sim \PP \left(X_A^{(i)}\,W_{h,i} > x\right)\,,
\eeam

From \eqref{eq.KLP.3.12} and the fact that the function $h_i$ is bounded from above (say by a constant $C_i \in (0,\,\infty)$, by similar reasons with \eqref{eq.KLP.2.A}), we obtain
\beam \label{eq.KLP.3.13}
\E\left[W_{h,i}^p \right] =\int_{s(G_{h,i})} y\,G_{h,i}(dy) =\int_{s(G_i)} y\,h_i(y)\,G_i(dy) \leq C_i\,\E\left[W_i^p \right] < \infty\,,
\eeam
where at the second step we used the fact that $s(G_{h,i}) = s(G_i)$, as it follows by 
relation \eqref{eq.KLP.3.6} and \eqref{eq.KLP.3.12}.

Since $X_A^{(i)} \stackrel{d}{\sim} F_A \in \mathcal{D}$, the $W_{h,i}$ is independent of $X_A^{(i)}$, and \eqref{eq.KLP.3.13} holds, for some $p> J_{F_A}^+$, by \cite[Th. 3.3(iv)]{cline:samorodnitsky:1994} we obtain
\beam \label{eq.KLP.3.14}
\PP\left(X_A^{(i)}\,W_{h,i} >x \right) \asymp \PP \left(X_A > x\right)\,.
\eeam
Hence, from relations \eqref{eq.KLP.3.11} and  \eqref{eq.KLP.3.14} we find
\beam \label{eq.KLP.3.15}
\PP\left(X_A^{(i)}\,W_{i} >x \right) \asymp \PP \left(X_A > x\right)\,.
\eeam
Further, due to \eqref{eq.KLP.3.13}, we have by \cite[Lem. 3.9]{tang:tsitsiashvili:2003} that $X_A^{(i)}\,W_{h,i} \stackrel{d}{\sim} H_{h,i}$, with upper Matuszewska index $J_{H_{h,i}}^+ =J_{F_A}^+$, while from \eqref{eq.KLP.3.11}, we see that $X_A^{(i)}\,W_{i} \stackrel{d}{\sim} H_A^{(i)}$, with upper Matuszewska index $J_{H_A^{(i)}}^+ =J_{H_{h,i}}^+$, that implies $J_{H_A^{(i)}}^+ =J_{F_A}^+$. Because of this, and since $W_{i-1}$ is independent of $X_A^{(i)}\,W_{i}$, with distribution $H_A^{(i)} \in \mathcal{D}$ (recall the characterization of $\mathcal{D}$ by the inequality $J_{H_{i}}^+ < \infty$), and again through \cite[Th. 3.3(iv)]{cline:samorodnitsky:1994}, taking into account the moment condition $\E\left[W_{i-1}^p \right] < \infty$ for some $p > J_{H_A^{(i)}}^+ = J_{F_A}^+$, through \eqref{eq.KLP.3.10}, we have
\beao
\PP\left(X_A^{(i)}\,W_{i}\,W_{i-1} >x \right) \asymp \PP \left(X_A > x\right)\,.
\eeao
Furthermore, again the upper Matuszewska index of $X_A^{(i)}\,W_{i}\,W_{i-1}$ coincides with $J_{H_A^{(i)}}^+ =J_{F_A}^+$, from \cite[Lem. 3.9]{tang:tsitsiashvili:2003}. 
Continuing on similar line by induction, we obtain relation \eqref{eq.KLP.3.10}, and consequently relation \eqref{eq.KLP.3.8} too.

From inclusion $F \in \mathcal{D}_A$ and since \eqref{eq.KLP.3.8} holds, we find that the distribution of the product ${\bf X}^{(i)}\,\Pi_i$ belongs to class $\mathcal{D}_A$, for any $i \in \bbn$ (see for example \cite[Prop. 2.2(ii)]{konstantinides:passalidis:2024g}). 

\item[(ii)] Relation \eqref{eq.KLP.3.B} can be directly proved through \eqref{eq.KLP.3.8}, since from \eqref{eq.KLP.3.C} we can apply immediately \cite[Lem. 4.5 (i)]{yang:gao:li:2016}.  

\item[(iii)]
From relation \eqref{eq.KLP.3.9} follows that it is enough to show that the sequence $\{X_A^{(i)}\,\Pi_i,\;i \in\bbn \}$ is  $TAI$. Without loss of generality, for any $i,\,j \in \bbn$, considering $i< j$, we obtain 
\beao
&&\lim_{x_i \wedge x_j \to \infty} \PP(X_A^{(i)} \,\Pi_i >x_i\,,\;X_A^{(j)} \,\Pi_j >x_j)\\[2mm] \notag
&&=\lim_{x_i \wedge x_j \to \infty}\PP(X_A^{(i)} \,\Pi_i >x_i\,,\;X_A^{(j)}\,W_j \,\Pi_{j-1} >x_j) =o\left[\PP\left(X_A^{(j)}\,\Pi_{j} >x_j \right) \right]\,.
\eeao
where at the second step we used \cite[Lem. 7]{tang:yuan:2014}, since from assertion (i) we have $H_A^{(j)} \in \mathcal{D}$, the fact that $X_A^{(j)}\,W_j$ is independent of $\Pi_{j-1}$ and of $X_A^{(i)} \,\Pi_i$ (since $i< j$), and also the moment condition 
\beao
\E\left[\Pi_{j-1}^p \right] = \prod_{k=1}^{j-1} \E\left[W_{k}^p \right] < \infty\,,
\eeao 
for some $p > J_{H_A^{(j)}}^+ = J_{F_A}^+$, as follows from the proof of assertion (i). Symmetrically, it is implied also the case $i >j$. So, we have shown that $\{X_A^{(i)}\,\Pi_i,\;i \in\bbn \}$ is  $TAI$.

Alternatively, the assertion (iii) is implied by similar steps of \cite[Lem. 4.6(i)]{yang:gao:li:2016}, without any obstacle from the fact that in that paper the $\{W_i\,,\;i \in \bbn\}$ should be identically distributed.    
\item[(iv)]
From assertion (i), since $F \in (\mathcal{D} \cap \mathcal{L})_A \subsetneq \mathcal{D}_A$.

From \eqref{eq.KLP.3.13} we have the moment condition $ \E\left[W_{h,i}^p \right] < \infty$, hence $X_A^{(i)}\,W_{h,i} \stackrel{d}{\sim} H_{h,i} \in \mathcal{D} \cap \mathcal{L} $ (see for example \cite[Cor. 5.2(c)]{leipus:siaulys:konstantinides:2023}). From \eqref{eq.KLP.3.11}, we obtain $H_A^{(i)} \in \mathcal{D} \cap \mathcal{L}$. Because of the equalities $J_{H_A^{(i)}}^+ =J_{H_{h,i}}^+ = J_{F_A}^+$, and the independence between $W_{i-1}$ and $X_A^{(i)}\,W_{i}$, by \cite[Cor. 5.2(c)]{leipus:siaulys:konstantinides:2023}, we obtain that the distribution of $X_A^{(i)}\,W_{i}\,W_{i-1}$ belongs to the class $\mathcal{D} \cap \mathcal{L}$, with upper Matuszewska index coinciding with $J_{F_A}^+$. By induction, we get that the distribution of  $X_A^{(i)}\,\Pi_i$ belongs to the class $\mathcal{D} \cap \mathcal{L}$, that through \eqref{eq.KLP.3.9} provides the desired result. ~\halmos
\end{enumerate}

The following lemma comes from \cite[Th. 4.1(ii)]{konstantinides:passalidis:2024g}.

\ble \label{lem.KLP.3.2}
Let $A \in \mathscr{R}$ be a fixed set. We suppose that the ${\bf Z}^{(1)},\,\ldots,\,{\bf Z}^{(n)}$ are non-negative, random vectors, that possess $TAI_A$, and their corresponding distributions $V_1,\,\ldots,\,V_n$ belong to the class $(\mathcal{D} \cap \mathcal{L})_A$. Then, it holds
\beao
 \PP\left(\sum_{i=1}^n {\bf Z}^{(i)} \in x\,A \right) \sim \sum_{i=1}^n \PP\left( {\bf Z}^{(i)} \in x\,A \right)\,. 
\eeao
\ele

Now we define the following randomly weighted sum 
\beam \label{eq.KLP.3.17}
{\bf S}_{n} =\sum_{i=1}^n {\bf X}^{(i)}\,\Pi_i\,,
\eeam
for any fixed $n \in \bbn$. The following lemma says that under the assumptions of Lemma \ref{lem.KLP.3.1}(iii), the sum from \eqref{eq.KLP.3.17} satisfies the multivariate linear single big jump principle. For similar statements, with several dependence structures, distribution classes and sets $A$, see \cite{resnick:willekens:1991}, \cite{chen:yang:2019}, \cite{cheng:konstantinides:wang:2024}, \cite{konstantinides:passalidis:2025o}, among others.

\ble \label{lem.KLP.3.3}
Let $A \in \mathscr{R}$ be a fixed set. We suppose that hold the assumptions of Lemma \ref{lem.KLP.3.1}, with $F \in (\mathcal{D} \cap \mathcal{L})_A$. Then for each $n \in \bbn$ it holds
\beam \label{eq.KLP.3.18}
\PP\left({\bf S}_{n} \in x\,A \right) \sim \sum_{i=1}^n \PP\left( {\bf X}^{(i)}\,\Pi_i \in x\,A \right) \sim \sum_{i=1}^n \PP\left( {\bf X}^{(i)}\,W_{h,i}\,\Pi_{i-1} \in x\,A \right)\,.
\eeam
\ele

\pr~
From Lemma \ref{lem.KLP.3.1}(iii), (iv), we obtain that the sequence $\{{\bf X}^{(i)}\,\Pi_i\,,\;i \in \bbn \}$ is $TAI_A$ and the distributions its terms belong to the class $(\mathcal{D} \cap \mathcal{L})_A$. Hence, the first relation in \eqref{eq.KLP.3.18} follows by application of Lemma \ref{lem.KLP.3.2}, with ${\bf Z}^{(i)}={\bf X}^{(i)}\,\Pi_i$, for any $i=1,\,\ldots,\,n$. The second relation is implied immediately through application of Lemma \ref{lem.KLP.3.1} (ii).
~\halmos

In the next lemma, we reformulate \cite[Lem. 3.3]{chen:yuan:2017}, see also \cite[Lem. 3.2]{hao:tang:2012}.

\ble \label{lem.KLP.3.4}
Let $Z$ be a real-valued random variable, with distribution $G \in \mathcal{D}\cap \mathcal{P_D}$. Then, for any pair $p_1,\,p_2$, such that $0 < p_1 < J_G^- \leq J_G^+ < p_2 < \infty$,  there exists a constant $C>0$ and some $x_0=x_0(p_1,\,p_2)>0$, such that for any non-negative random variable $\xi$, independent of $Z$, it holds
\beao
\dfrac{\PP\left(\xi\,Z > x \right)}{\PP\left(Z > x \right)} \leq C\,\E\left[ \xi^{p_1} \vee \xi^{p_2} \right]\,,
\eeao
for any $x > x_0$.
\ele

In the following lemma we find out that relation \eqref{eq.KLP.3.18} holds uniformly with respect to $n \in \bbn$, under some stricter conditions, in comparison to these of Lemma \ref{lem.KLP.3.3}. Namely, we show that
\beam \label{eq.KLP.3.20} \notag
&&\lim_{\xto} \sup_{n \in \bbn} \left| \dfrac{\PP\left({\bf S}_{n} \in x\,A\right)}{\sum_{i=1}^n \PP\left( {\bf X}^{(i)}\,\Pi_i \in x\,A \right)} -1 \right|=\lim_{\xto} \sup_{n \in \bbn} \left| \dfrac{\PP\left({\bf S}_{n} \in x\,A\right)}{\sum_{i=1}^n \PP\left( {\bf X}^{(i)}\,W_{h,i}\,\Pi_{i-1} \in x\,A \right)} -1 \right| \\
&&=0\,.
\eeam
Here, for the uniformity of the convergence, with respect to $n \in \bbn$, we require the inclusion $F \in (\mathcal{D} \cap \mathcal{A})_A$, instead of $F \in (\mathcal{D} \cap \mathcal{L})_A$ from Lemma \ref{lem.KLP.3.3}, and a more strict moment condition on $W_i$, that depends also on two Matuszewska indexes of $F_A$. The restriction from class $(\mathcal{D} \cap \mathcal{L})_A$ to $(\mathcal{D} \cap \mathcal{A})_A$ seems small, but it is in fact also necessary, in order to get for ${\bf S}_{\infty}$ a non-defective distribution. This can be seen from relation \eqref{eq.KLP.3.24} in the next lemma.

\ble \label{lem.KLP.3.5}
Let $A \in \mathscr{R}$ be a fixed set. We suppose that Assumption \ref{ass.KLP.3.1}, is satisfied for $\{({\bf X}^{(i)},\,W_i)\,,\;i \in \bbn\}$, with $F \in (\mathcal{D} \cap \mathcal{A})_A$. If there exists a pair $p_1,\,p_2$, with $0 < p_1 < J_{F_A}^- \leq J_{F_A}^+ < p_2 < \infty$, such that it holds $\E\left[ W_i^{p_1} \vee W_i^{p_2} \right]< 1$, for any $i \in \bbn$ then relation \eqref{eq.KLP.3.20} is true.
\ele

\pr~
Let $M \in \bbn$. Then from Lemma \ref{lem.KLP.3.3}, we obtain that \eqref{eq.KLP.3.20} is true for any $n \leq M$. We proceed now to show that it holds also for any $n > M$. We prove only that the first relation in \eqref{eq.KLP.3.20} equals to zero, since after this the proof of the second relation follows directly by Lemma \ref{lem.KLP.3.1} (ii). 

Further, we need to prove two intermediate auxiliary relations. The first one is
\beam \label{eq.KLP.3.21}
\lim_{M \to \infty} \limsup \dfrac{\sum_{i=M+1}^{\infty} \PP\left( {\bf X}^{(i)}\,\Pi_i \in x\,A \right)}{\PP\left( {\bf X}^{(1)}\,W_1 \in x\,A \right)} =0\,.
\eeam
Initially, by \cite[Cor. 5.2(c)]{leipus:siaulys:konstantinides:2023} and \cite[Th. 5.1(i)]{konstantinides:passalidis:2025b}, we have $X_A^{(i)}\,W_{h,i} \stackrel{d}{\sim} H_{h,i} \in \mathcal{D}\cap \mathcal{A}$, and by relation  \eqref{eq.KLP.3.11}, we get $X_A^{(i)}\,W_{i} \stackrel{d}{\sim} H_A^{(i)} \in \mathcal{D}\cap \mathcal{A} \subsetneq \mathcal{D}\cap \mathcal{P_D}$.  As consequence, by Lemma \ref{lem.KLP.3.4}, for any $i \in \bbn$ there exist $C \in (0,\,\infty)$ and $x_0>0$, such that it holds
\beam \label{eq.KLP.3.22} \notag
&&\dfrac{\PP\left( {\bf X}^{(i)}\,\Pi_i \in x\,A \right)}{\PP\left( {\bf X}^{(1)}\,W_1 \in x\,A \right)} = \dfrac{\PP\left( X_A^{(i)}\,W_i\,\Pi_{i-1} > x \right)}{\PP\left(X_A^{(1)}\,W_1 > x \right)} \leq C\,\E\left[ \Pi_{i-1}^{p_1} \vee \Pi_{i-1}^{p_2} \right] \, \dfrac{\PP\left( X_A^{(i)}\,W_i > x \right)}{\PP\left(X_A^{(1)}\,W_1 > x \right)} \\[2mm] 
&& \leq C\,C^* \,\E\left[ \Pi_{i-1}^{p_1} \vee \Pi_{i-1}^{p_2} \right] \leq C\,C^* \,\left( \E\left[ \Pi_{i-1}^{p_1}\right] + \E\left[ \Pi_{i-1}^{p_2} \right]\right)\,,
\eeam
for any $x> x_0$, where the constant $C^*>0$, comes from Lemma \ref{lem.KLP.3.1}(i), see relation \eqref{eq.KLP.3.15}. We denote
\beao
\E\left[\check{W}_i^{p_k} \right]:=\bigvee_{j=1}^{i-1} \E\left[W_{j}^{p_k} \right]\,,
\eeao
for $k=1,\,2$. Then, since the $\{W_i\,,\;i \in \bbn\}$ are independent each other, we obtain
\beam \label{eq.KLP.3.23}
\E\left[\Pi_{i-1}^{p_k} \right]=\prod_{j=1}^{i-1} \E\left[W_{j}^{p_k} \right] \leq \left( \E\left[\check{W}_i^{p_k} \right]\right)^{i-1}\,,
\eeam
for $k=1,\,2$. From \eqref{eq.KLP.3.22} and  \eqref{eq.KLP.3.23}, taking into consideration that $\E\left[ W_i^{p_1} \vee W_i^{p_2} \right]< 1$, for all $i \in \bbn$, then for any $\vep>0$ we can find large enough $M=M(\vep) \in \bbn$, such that it holds
\beao
\dfrac{\sum_{i=M+1}^{\infty} \PP\left( {\bf X}^{(i)}\,\Pi_i \in x\,A \right)}{\PP\left( {\bf X}^{(1)}\,W_1 \in x\,A \right)} < \vep\,.
\eeao  
Hence, by the arbitrariness in the choice of $\vep>0$, we obtain that \eqref{eq.KLP.3.21} is true. 

The second auxiliary relation, which we shall prove is the following
\beam \label{eq.KLP.3.24}
\PP\left(\sum_{i=1}^{\infty} {\bf X}^{(i)}\,\Pi_i \in x\,A \right)= O\left[\PP\left( {\bf X}\,W_1 \in x\,A \right) \right] \,.
\eeam
From the moment condition for the $\{W_i\,,\;i \in \bbn\}$, we can find some $\delta \in (0,\,1)$, arbitrarily close to unity, such that it holds
\beam \label{eq.KLP.3.25}
\dfrac {\E\left[W_i^{p_1} \right]}{\delta^{p_1}} \bigvee \dfrac {\E\left[W_i^{p_2} \right]}{\delta^{p_2}} <1\,.
\eeam
Then we obtain
\beam \label{eq.KLP.3.26} \notag
&&\PP\left(\sum_{i=1}^{\infty} {\bf X}^{(i)}\,\Pi_i \in x\,A \right)= \PP\left(\sup_{{\bf p}\in I_A} {\bf p}^T\,\left(\sum_{i=1}^{\infty} {\bf X}^{(i)}\,\Pi_i \right) > x \right)\\[2mm]
&& \leq \PP\left(\sum_{i=1}^{\infty}  X_{A}^{(i)}\,\Pi_i >x \right)=  \PP\left(\sum_{i=1}^{\infty}  X_{A}^{(i)}\,\Pi_i >(1-\delta)\,x\,\sum_{i=1}^{\infty} \delta^{i-1} \right)\\[2mm] \notag
&& \leq \sum_{i=1}^{\infty} \PP\left( X_{A}^{(i)}\,\Pi_i >(1-\delta)\,\delta^{i-1}\,x \right)= \sum_{i=1}^{\infty}  \PP\left( X_{A}^{(i)}\,W_i\,\prod_{j=1}^{i-1}\dfrac{W_j}{\delta} >(1-\delta)\,x \right)\,,
\eeam
where at the second step we used the monotonicity of the event sequences 
\beao
\left[ \left\{\sup_{{\bf p} \in I_A} {\bf p}^{\top}\left(\sum_{i=1}^n {\bf X}^{(i)}\,\Pi_i\right) > x\right\} \,,\;n \in \bbn\right]\,, \qquad \left[ \left\{\sum_{i=1}^n  X_A^{(i)}\,\Pi_i > x \right\}\,,\;n \in \bbn\right]\,.
\eeao

Because of the inclusion $H_{A}^{(i)}\in \mathcal{D} \cap \mathcal{A}$, we apply to each term of the last sum in \eqref{eq.KLP.3.26} Lemma \ref{lem.KLP.3.4}, separately, as before with relation \eqref{eq.KLP.3.22}, and we have that for any large enough $x>0$ the inequalities
\beam \label{eq.KLP.3.27} \notag
&&\PP\left( X_{A}^{(i)}\,W_i\,\prod_{j=1}^{i-1}\dfrac{W_j}{\delta} >(1-\delta)\,x \right)\\[2mm]
&&\leq C\,\E\left[\left(\prod_{j=1}^{i-1} \dfrac{W_j^{p_1} }{\delta}\right)^{p_1} \bigvee \left(\prod_{j=1}^{i-1} \dfrac {W_i^{p_2} }{\delta}\right)^{p_2}\right]\,\PP\left( X_{A}^{(i)}\,W_i >(1-\delta)\,x \right) \\[2mm] \notag
&&\leq C\,C^*\,\PP\left( X_{A}^{(1)}\,W_1 >x \right)\,\left[\left(\dfrac{\E\left[\check{W}_i^{p_1} \right]}{\delta^{p_1}}\right)^{i-1}+ \left(\dfrac{\E\left[\check{W}_i^{p_2} \right]}{\delta^{p_2}}\right)^{i-1}\right]\,,
\eeam
hold, where the constant $C^*>0$ stems from Lemma \ref{lem.KLP.3.1}(i), relation \eqref{eq.KLP.3.15}, and the inclusion $H_{A}^{(i)}\in \mathcal{D}$. Therefore, by relations \eqref{eq.KLP.3.26} and \eqref{eq.KLP.3.27}, we find
\beao
\PP\left(\sum_{i=1}^{\infty} {\bf X}^{(i)} \Pi_i \in x\,A \right) \leq C C^* \PP\left( X_{A}^{(1)} W_1 >x \right) \sum_{i=1}^{\infty} \left[\left(\dfrac{\E\left[\check{W}_i^{p_1} \right]}{\delta^{p_1}}\right)^{i-1}+ \left(\dfrac{\E\left[\check{W}_i^{p_2} \right]}{\delta^{p_2}}\right)^{i-1}\right] 
\eeao
with the last sum bounded from above by the quantity
\beao
2\,\left[\left(1-\dfrac{\E\left[\check{W}_i^{p_1} \right]}{\delta^{p_1}}\right)^{-1} \bigvee \left(1-\dfrac{\E\left[\check{W}_i^{p_2} \right]}{\delta^{p_2}}\right)^{-1}\right] < \infty\,,
\eeao
(recall relation  \eqref{eq.KLP.3.25}). From these two last relation we obtain  \eqref{eq.KLP.3.24}.

Let us introduce a new random variable $\eta$, that is non-negative, independent of the other sources of randomness, and such that it holds
\beam \label{eq.KLP.3.28} \notag
&& \PP(\eta >x) \sim c\,\PP\left( X_{A}^{(1)}\,W_1> x \right)\\[2mm]
&&\PP\left(\sum_{i=1}^{\infty} X_{A}^{(i)}\,\Pi_i> y \right) \leq \PP\left(\eta >y\right)\,,
\eeam
for some constant $c>0$, and for any $y \in \bbr$, respectively. We notice that the second relation in \eqref{eq.KLP.3.28} is feasible, due to relation  \eqref{eq.KLP.3.24}. Hence, for any $n> M$ it holds
\beam \label{eq.KLP.3.29} \notag
&&\PP\left({\bf S}_{n} \in x\,A\right) \leq \PP\left(\sum_{i=1}^{n} X_{A}^{(i)}\,\Pi_i> x \right) \notag\\[2mm] 
&&= \PP\left(\sum_{i=1}^{M}  X_{A}^{(i)}\,\Pi_i + \left(\sum_{i=M+1}^{n} X_{A}^{(i)}\,\prod_{j=M+1}^i W_i \right)\,\Pi_M> x \right) \\[2mm] \notag
&&\leq \PP\left(\sum_{i=1}^{M}  X_{A}^{(i)}\,\Pi_i + \eta\,\Pi_M> x \right) \sim \sum_{i=1}^{M}\PP\left(  X_{A}^{(i)}\,\Pi_i > x \right) + \PP\left(\eta\,\Pi_M> x \right)\,,
\eeam
where at the first step we used \cite[Prop. 2.4]{konstantinides:passalidis:2024g}, at the third step we take into account the second relation of  \eqref{eq.KLP.3.28}, while at the last step we use the $TAI$ dependence of the $\{ X_{A}^{(1)}\,\Pi_1,\,\ldots,\, X_{A}^{(M)}\,\Pi_M,\,\eta\,\Pi_M\}$, with distributions from class $\mathcal{D} \cap \mathcal{L}$. Indeed, the  $TAI$ dependence of $\{ X_{A}^{(1)}\,\Pi_1,\,\ldots,\, X_{A}^{(M)}\,\Pi_M\}$ is implied by Lemma \ref{lem.KLP.3.1} (iii). So we have to show that for any $i = 1,\,\ldots,\,M$ the $X_{A}^{(i)}\,\Pi_i$ and $\eta\,\Pi_M$ are $TAI$. Since $\eta$ is independent of any other source of randomness, and by the first relation of \eqref{eq.KLP.3.28}, the distribution of $\eta$ belongs to class $\mathcal{D}$, and the upper Matuszewska index $J_{\eta}^+$ coincides with the $J_{H_A^{(i)}}^+ = J_{F_A}^+$, while $\E\left[ \Pi_M^p\right] < \infty$, for some $p> J_{\eta}^+ =J_{F_A}^+$, and hence we can apply \cite[Lem. 7]{tang:yuan:2014}, to obtain
\beam \label{eq.KLP.3.29a}
\lim_{x_i \wedge x_j \to \infty} \dfrac{\PP\left( X_{A}^{(i)}\,\Pi_i > x_i\,,\;\eta\,\Pi_M > x_j\right)}{\PP\left( \eta\,\Pi_M > x_j\right)} =0\,.
\eeam 
In fact the limit in \eqref{eq.KLP.3.29a} holds also for $x=x_i=x_j$ (which represents a special case of $x_i \wedge x_j$, that is required in $TAI$).
The symmetric limit 
\beao
\lim_{x \to \infty} \dfrac{\PP\left( X_{A}^{(i)}\,\Pi_i > x\,,\;\eta\,\Pi_M > x\right)}{\PP\left( X_{A}^{(i)}\,\Pi_i > x\right)} \,,
\eeao
equals also to zero, due to relation \eqref{eq.KLP.3.29a} in combination with 
\beao
\PP(\eta\,\Pi_M >x) \asymp \PP(\eta >x)\,,
\eeao 
from \cite[Th. 3.3(iv)]{cline:samorodnitsky:1994}, the first relation in \eqref{eq.KLP.3.28} and relation \eqref{eq.KLP.3.8}. That means, the $\{ X_{A}^{(1)}\,\Pi_1,\,\ldots,\, X_{A}^{(M)}\,\Pi_M,\,\eta\,\Pi_M\}$ possess the  $TAI$ dependence, and from Lemma \ref{lem.KLP.3.1} (iv), together with the first relation  of \eqref{eq.KLP.3.28} their distribution belongs to class  $\mathcal{D} \cap \mathcal{L}$. Therefore, at the last step of \eqref{eq.KLP.3.29} we can apply \cite[Th. 3.1]{geluk:tang:2009}.

Further, since from the first relation  of \eqref{eq.KLP.3.28} follows that the distribution of $\eta$ belongs to $\mathcal{D} \cap \mathcal{A}$, hence by Lemma \ref{lem.KLP.3.4} we obtain that for any $\vep>0$, there exists some $M=M(\vep) \in \bbn$ such that it holds
\beam \label{eq.KLP.3.30}
&&\dfrac{\PP\left(\eta\,\Pi_M> x \right)}{\PP\left( X_{A}^{(1)}\,W_1 > x \right)} \sim \dfrac{\PP\left(\eta\,\Pi_M> x \right)}{\PP\left( \eta > x \right)/c} \leq c\,C\, \E\left[ \Pi_{M}^{p_1}\vee \Pi_{M}^{p_2} \right] \\[2mm]  \notag
&&\leq c\,C\,\left( \prod_{j=1}^M \left\{\E\left[ W_{j}^{p_1}\right] + \E\left[ W_{j}^{p_2} \right] \right\} \right) \leq c\,C\,\left[\left( \E\left[ \check{W}_{M+1}^{p_1}\right]\right)^M +\left( \E\left[ \check{W}_{M+1}^{p_2} \right]\right)^M \right] \leq \vep\,.
\eeam
From relations \eqref{eq.KLP.3.29} and \eqref{eq.KLP.3.30}, and for all $n>M$ we find
\beam \label{eq.KLP.3.31} \notag
&&\PP\left({\bf S}_{n} \in x\,A\right) \lesssim \sum_{i=1}^M \PP\left(X_{A}^{(i)}\,\Pi_i > x \right) + \vep\,\PP\left(X_{A}^{(1)}\,W_1 > x \right)\\[2mm] 
&& \leq (1+\vep)\, \sum_{i=1}^M \PP\left(X_{A}^{(i)}\,\Pi_i > x \right)\leq (1+\vep)\, \sum_{i=1}^n \PP\left({\bf X}^{(i)}\,\Pi_i \in x\,A\right)\,.
\eeam

From the other hand side, since set $x\,A$ is increasing, and the summands of the ${\bf S}_{n}$ are non-negative, we obtain, for all $n>M$
\beam \label{eq.KLP.3.32}
&&\PP\left({\bf S}_{n} \in x\,A\right) \geq \PP\left({\bf S}_{M} \in x\,A\right) \sim \left( \sum_{i=1}^{n} - \sum_{i=M+1}^n \right)\PP\left({\bf X}^{(i)}\,\Pi_i \in x\,A\right) \\[2mm]  \notag
&& \geq \sum_{i=1}^n \PP\left({\bf X}^{(i)}\,\Pi_i \in x\,A\right) - \vep\,\PP\left({\bf X}^{(1)}\,W_1 \in x\,A\right) \geq (1-\vep)\,\sum_{i=1}^n \PP\left({\bf X}^{(i)}\,\Pi_i \in x\,A\right)\,,
\eeam
 where at the second we apply Lemma \ref{lem.KLP.3.3}, and at the third step we use relation \eqref{eq.KLP.3.21}.  

From relations \eqref{eq.KLP.3.31}, \eqref{eq.KLP.3.32} and the arbitrary choice of $\vep>0$, follows that the first relation of \eqref{eq.KLP.3.20} tends to zero for all $n> M$, namely we have that  \eqref{eq.KLP.3.18} holds uniformly for $n>M$.
~\halmos

The following lemma plays crucial role for the proof of Corollary \ref{cor.KLP.3.a} and in fact provides more explicit expressions than \eqref{eq.KLP.3.20} that was proved in Lemma \ref{lem.KLP.3.5}, when we restrict $F$ in class $MRV$, under some slightly weaker moment conditions on $\{W_i\,,\; i \in \bbn\}$.

\ble \label{lem.KLP.3.6}
Let $A \in \mathscr{R}$ some fixed set. Let assume that Assumption \ref{ass.KLP.3.1} is satisfied for $\{({\bf X}^{(i)},\,W_i)\,,\;i \in \bbn\}$, with $F \in MRV(\alpha,\,\mu)$, and $\alpha \in (0,\,\infty)$. If for any $i \in \bbn$ the inequalities
\beam \label{eq.KLP.3.c}
\E[W_i^{\alpha}] < 1\,, \qquad \E[W_i^{p_2}] < \infty\,, 
\eeam
are true, for some $p_2 > \alpha$, then it holds
\beam \label{eq.KLP.3.d}
\lim_{\xto} \sup_{n \in \bbn} \left|\dfrac{\PP[{\bf S}_n \in x\,A]}{\mu(A)\,\bV(x)\,\sum_{i=1}^n \E[W_{h,i}^{\alpha}]\,\E[\Pi_{i-1}^{\alpha}] } - 1 \right|= 0\,. 
\eeam
\ele

\pr~
Firstly, we note that since $F \in MRV(\alpha,\,\mu)$, then for any $A \in \mathscr{R}$, it holds $F_A \in \mathcal{R}_{-\alpha}$ (see proof of \cite[Prop. 4.14]{samorodnitsky:sun:2016}), and therefore we obtain $J_{F_A}^- = J_{F_A}^+ =\alpha$. 

From \eqref{eq.KLP.3.c}, via monotone convergence theorem we find
\beao
&&\lim_{\gamma \to 0} \E[W_i^{\alpha - \gamma} \vee W_i^{\alpha + \gamma}]= \lim_{\gamma \to 0} \int_0^{\infty} (y^{\alpha - \gamma}  \vee y^{\alpha + \gamma})\,\PP(W_i \in dy) \\[2mm]
&& = \int_0^{\infty} \lim_{\gamma \to 0}(y^{\alpha - \gamma}  \vee y^{\alpha + \gamma})\,\PP(W_i \in dy) = \E[W_i^{\alpha}] < 1\,.
\eeao
Hence, there exists some $\delta >0$ such that $\E[W_i^{\alpha - \delta} \vee W_i^{\alpha + \delta}] < 1$, and so by Lemma \ref{lem.KLP.3.5} follows that relation \eqref{eq.KLP.3.20} is valid.

For any $i \in \bbn$, from the second relation in \eqref{eq.KLP.3.c} and relation \eqref{eq.KLP.2.A}, via the application of Breiman's theorem (see, \cite[Prop. 5.2 (iv)]{leipus:siaulys:konstantinides:2023}) we obtain
\beam \label{eq.KLP.3.e} \notag
&&\PP({\bf X}^{(i)}\,\Pi_i \in x\,A) \sim \PP({\bf X}^{(i)}\,W_{h,i}\,\Pi_{i-1} \in x\,A) = \PP( X_A^{(i)}\,W_{h,i}\,\Pi_{i-1} > x) \\[2mm]
&&\sim  \PP( X_A^{(i)}> x) \E[W_{h,i}^{\alpha}]\,\E[\Pi_{i-1}^{\alpha}] \sim \mu(A)\,\bV(x)\,\E[W_{h,i}^{\alpha}]\,\E[\Pi_{i-1}^{\alpha}]  \,, 
\eeam
where at the last step we used that $F \in MRV(\alpha,\,\mu)$.

From \eqref{eq.KLP.3.e} and \eqref{eq.KLP.3.18}, we find that for any fixed $n \in \bbn$ it holds
\beam \label{eq.KLP.3.f} \notag
\PP({\bf S}_n \in x\,A) \sim \mu(A)\,\bV(x)\,\sum_{i=1}^n \E[W_{h,i}^{\alpha}]\,\E[\Pi_{i-1}^{\alpha}]  \,. 
\eeam 

Now it remains to show that \eqref{eq.KLP.3.d} holds. From the proof of Lemma \ref{lem.KLP.3.5} (recall \eqref{eq.KLP.3.31}, \eqref{eq.KLP.3.32}), for any $\vep \in (0,\,1)$ there exists some large enough $M \in\bbn$, such that for any $n > M$ it holds
\beao
&&(1-\vep)\, \sum_{i=1}^M \PP({\bf X}^{(i)}\,\Pi_i \in x\,A) \lesssim \PP({\bf S}_n \in x\,A) \\[2mm]
&&\lesssim (1+\vep)\,\sum_{i=1}^M \PP({\bf X}^{(i)}\,\Pi_i \in x\,A) \,.
\eeao
From the last relation, in combination with \eqref{eq.KLP.3.e} we get
\beam \label{eq.KLP.3.g} \notag
&&(1-\vep)\,\mu(A)\,\bV(x)\,\sum_{i=1}^M \E[W_{h,i}^{\alpha}]\,\E[\Pi_{i-1}^{\alpha}]\lesssim \PP({\bf S}_n \in x\,A) \\[2mm]
&&\lesssim (1+\vep)\,\mu(A)\,\bV(x)\,\sum_{i=1}^M \E[W_{h,i}^{\alpha}]\,\E[\Pi_{i-1}^{\alpha}] \,. 
\eeam
From \eqref{eq.KLP.3.g} we immediately obtain that for all $n > M$ from right hand side, it holds
\beam \label{eq.KLP.3.h}
\PP({\bf S}_n \in x\,A) \lesssim (1+\vep)\,\mu(A)\,\bV(x)\,\sum_{i=1}^n \E[W_{h,i}^{\alpha}]\,\E[\Pi_{i-1}^{\alpha}] \,. 
\eeam
From the left hand side, due to \eqref{eq.KLP.2.A}, \eqref{eq.KLP.3.c} and the fact that the $\{W_i\,,\;i \in\bbn\}$ are independent, we find that the series $\sum_{i=1}^{\infty} \E[W_{h,i}^{\alpha}]\,\E[\Pi_{i-1}^{\alpha}]$ is convergent, hence by \eqref{eq.KLP.3.g} we conclude that for all $n > M$ it holds
\beam \label{eq.KLP.3.i} \notag
&&\PP({\bf S}_n \in x\,A) \gtrsim (1-\vep)\,\mu(A)\,\bV(x)\,\left(\sum_{i=1}^n - \sum_{i=M+1}^{\infty} \right) \E[W_{h,i}^{\alpha}]\,\E[\Pi_{i-1}^{\alpha}] \\[2mm]
&&\gtrsim (1-\vep)^2\,\mu(A)\,\bV(x)\,\sum_{i=1}^n \E[W_{h,i}^{\alpha}]\,\E[\Pi_{i-1}^{\alpha}] \,. 
\eeam
From relations \eqref{eq.KLP.3.h}, \eqref{eq.KLP.3.i} and the arbitrariness in the choice of $\vep > 0$, we obtain that 
\beao
\lim_{\xto} \sup_{N<n \in \bbn} \left| \dfrac{\PP({\bf S}_n \in x\,A)}{\mu(A)\,\bV(x)\sum_{i=1}^n \E[W_{h,i}^{\alpha}]\,\E[\Pi_{i-1}^{\alpha}] } -1\right| = 0\,.
\eeao
Further, from \eqref{eq.KLP.3.18} and \eqref{eq.KLP.3.e} follows immediately
\beao
\lim_{\xto} \sup_{N\geq n \in \bbn} \left| \dfrac{\PP({\bf S}_n \in x\,A)}{\mu(A)\,\bV(x)\sum_{i=1}^n \E[W_{h,i}^{\alpha}]\,\E[\Pi_{i-1}^{\alpha}] } -1\right| = 0\,.
\eeao
From these two last relations we obtain that \eqref{eq.KLP.3.d} is valid.
~\halmos

Now, we are ready to prove Theorem \ref{th.KLP.3.1}.

\noindent{\bf Proof of Theorem \ref{th.KLP.3.1}}~
For any $i \in \bbn$ we define 
\beao
W_i = e^{\xi(\tau_{i-1})-\xi(\tau{i})}\,.
\eeao 
Firstly, from Assumption \ref{ass.KLP.2.1} the $\{W_i \,,\;i \in \bbn\}$ are independent each other. Further, by relation \eqref{eq.KLP.3.1} there exists a pair $p_1,\,p_2$, with $0< p_1 < J_{F_A}^- \leq J_{F_A}^+ < p_2 < \infty$ such that the inequality 
\beao
\E\left[W_i^{p_1} \vee W_i^{p_2} \right]< 1\,,
\eeao 
is satisfied for any $i \in \bbn$. Hence, for any $i \in \bbn$, it holds
\beam \label{eq.KLP.3.33}
\PP\left({\bf X}^{(i)} \in x\,A\;\big|\; e^{\xi(\tau_{i-1}) - \xi(\tau_i)}=y\right) \sim h_i(y)\,\PP\left({\bf X}^{(i)} \in x\,A\right)\,,
\eeam 
uniformly with respect to $y \in s(Q_i)$, with function $h_i$ to satisfy relation \eqref{eq.KLP.3.6}. Further, by Assumption \ref{ass.KLP.2.2}, due to relation \eqref{eq.KLP.2.9}, we get that \eqref{eq.KLP.3.6} is true. We also denote the sequence of independent (from other sources of randomness) $\{W_{h,i}\,,\; i \in \bbn\}$ by
\beao
W_{h,i} = e^{\xi_{h_i}(\tau_{i-1}^*)-\xi_{h_i}(\tau{i}^*)}\,.
\eeao

Furthermore, since
\beao
\Pi_i:= e^{ -\xi(\tau_i)}=\prod_{j=1}^i e^{\xi(\tau_{j-1}) - \xi(\tau_j)}\,,
\eeao 
and because of \eqref{eq.KLP.3.33} and the moment condition for the $\{W_i \,,\;i \in \bbn\}$, we can apply Lemma \ref{lem.KLP.3.5} with $n=\infty$ to find that
\beao
&&\PP\left({\bf D}(\infty) \in x\,A \right) = \PP\left(\sum_{i=1}^{\infty} {\bf X}^{(i)}\,e^{-\xi(\tau_i)} \in x\,A \right) \sim \sum_{i=1}^{\infty}\PP\left( {\bf X}^{(i)}\,e^{-\xi(\tau_i)} \in x\,A \right) \\[2mm]
&& \sim \sum_{i=1}^{\infty} \PP\left( {\bf X}^{(i)}\,e^{\xi_{h_i}(\tau_{i-1}^*)-\xi_{h_i}(\tau_{i}^*)}\,e^{-\xi(\tau_{i-1})} \in x\,A \right)\,,
\eeao 
which provides relation \eqref{eq.KLP.3.2}.
~\halmos

\noindent{\bf Proof of Corollary \ref{cor.KLP.3.a}}~
It follows by similar way as the proof of Theorem \ref{th.KLP.3.1}, with the only difference that we apply Lemma \ref{lem.KLP.3.6} instead of Lemma \ref{lem.KLP.3.5}.
~\halmos

\noindent{\bf Proof of Corollary \ref{cor.KLP.3.1}}~
\begin{enumerate}
\item[(i)]
Because of relation  \eqref{eq.KLP.3.2}, and since the process $\{N(t)\,,\;t\geq 0\}$ is independent of any other source of randomness, we obtain by Theorem \ref{th.KLP.3.1}
\beao
&&\PP\left({\bf D}(\infty) \in x\,A \right) =\sum_{i=1}^{\infty} \PP\left({\bf X}^{(i)}\,e^{-\xi(\tau_i)} \in x\,A \right) = \sum_{i=1}^{\infty}\PP\left( {\bf X}^{(i)}\,e^{-\xi(\tau_i)} \in x\,A\,,\;\tau_i < \infty \right)\\[2mm]
&&=\int_0^{\infty} \PP\left({\bf X}^{(i)}\,e^{-\xi(s)} \in x\,A \right)\,\lambda(ds) \,,
\eeao
which provides the first relation in \eqref{eq.KLP.3.5}. 

For the second relation in \eqref{eq.KLP.3.5} we obtain
\beam \label{eq.KLP.3.bb} \notag
&&\PP\left({\bf D}(\infty) \in x\,A \right) \sim \sum_{i=1}^{\infty} \PP\left({\bf X}^{(i)}\,e^{\xi_{h_i}(\tau_{i-1}^*)-\xi_{h_i}(\tau_{i}^*)}\,e^{-\xi(\tau_{i-1})} \in x\,A \right) \\[2mm] \notag
&&= \sum_{i=1}^{\infty} \int_{K} \int_{\Lambda} \PP\left( {\bf X}^{(i)}\,e^{\xi_{h_i}(s_{i-1})-\xi_{h_i}(s_{i})}\,e^{-\xi(t)} \in x\,A \right)\,\PP(\tau_{i-1}^* \in ds_{i-1}\,,\;\tau_{i}^* \in ds_{i})\,\PP(\tau_{i-1}\in dt)\\[2mm]
&&=\sum_{i=1}^{\infty} \int_{K} \PP\left({\bf X}^{(i)}\,e^{\xi_{h_i}(s_{i-1})-\xi_{h_i}(s_{i})}\,e^{-\xi(t)} \in x\,A \right)\,\PP(\tau_{i-1}^* \in ds_{i-1}\,,\;\tau_{i}^* \in ds_{i})\,\lambda(dt)\,.
\eeam 
\item[(ii)]
Relation \eqref{eq.KLP.3.a} follows through Corollary \ref{cor.KLP.3.a}, via a similar way as with \eqref{eq.KLP.3.bb}.
\end{enumerate}
~\halmos

\noindent \textbf{Acknowledgments.} 
We feel the pleasant duty to express our sincere gratitude to Prof. Jinzhu Li for his useful comments that improved significantly the text.

\noindent \textbf{Funding Declaration.} 
No competing or funding interests, that influence the results of this paper.

\end{document}